\documentclass[12pt]{amsart}
\usepackage{amscd, amsfonts, amssymb}
\usepackage{amscd}
\usepackage{epsfig}
\usepackage{amssymb}
\usepackage[mathscr]{eucal}
\usepackage{verbatim}
\usepackage{fullpage}
\usepackage{latexsym}
\usepackage{lscape}

\newcommand\bb{\mathfrak b}
\newcommand\cc{\mathfrak c}

\renewcommand\gg{\mathfrak g}
\newcommand\hh{\mathfrak h}
\newcommand\kk{\mathfrak k}
\newcommand\frakl{\mathfrak l}
\newcommand\mm{\mathfrak m}
\newcommand\nn{\mathfrak n}
\newcommand\pp{\mathfrak p}

\newcommand\frakt{\mathfrak t}
\newcommand\zz{\mathfrak z}
\newcommand\mf{\mathfrak}

\newcommand\uu{\mathfrak u}

\newcommand{\inprod}[1]{\langle #1\rangle}

\newcommand\inverse{{^{-1}}}
\renewcommand{\check}{^{\vee}}

\newcommand\gl{\mathfrak{gl}}
\newcommand\Sl{\mathfrak{sl}}
\newcommand\ra{\rightarrow}

\DeclareMathOperator{\Ad}{Ad}

\DeclareMathOperator{\Char}{char}

\DeclareMathOperator{\GL}{GL}

\DeclareMathOperator{\SL}{SL}
\DeclareMathOperator{\PGL}{PGL}

\DeclareMathOperator{\SO}{SO}

\DeclareMathOperator{\SP}{Sp}

\DeclareMathOperator{\Hom}{Hom}
\DeclareMathOperator{\Lie}{Lie}

\numberwithin{equation}{section}

\newtheorem{thm}[equation]{Theorem} 
 \newtheorem{lem}[equation]{Lemma}
 \newtheorem{cor}[equation]{Corollary}
\newtheorem{prop}[equation]{Proposition}

 \theoremstyle{definition}
\newtheorem{defn}[equation]{Definition}

\newtheorem{exmp}[equation]{Example} \theoremstyle{remark}
\newtheorem{rem}[equation]{Remark} 
\theoremstyle{remark} \newtheorem{rems}[equation]{Remarks}
\newtheorem{question}[equation]{Question}

\thanks{2000 {\it Mathematics Subject Classification}.
20G15, 14L24.}
\keywords{$G$-complete reducibility, separability, reductive pairs}

\title[Complete Reducibility and Separability] {Complete Reducibility
and Separability}
\author[M.\ Bate]{Michael Bate}
\address%[M.\  Bate]
{Christ Church College, Oxford University, Oxford, OX1 1DP, United Kingdom}
\email{bate@maths.ox.ac.uk}

\author[B.\ Martin]{Benjamin Martin}
\address%[B.\ Martin]
{Mathematics and Statistics Department, University of Canterbury,
Private Bag 4800, Christchurch 1, New Zealand}
\email{B.Martin@math.canterbury.ac.nz}

\author[G.\ R\"ohrle]{Gerhard R\"ohrle}
\address%[G.~R\"{o}hrle]
{Fakult\"at f\"ur Mathematik,
Ruhr-Universit\"at Bochum,
D-44780 Bochum, Germany}
\email{gerhard.roehrle@rub.de}

\author[R.\ Tange]{Rudolf Tange}
\address%[R.~Tange]
{Fakult\"at f\"ur Mathematik,
Ruhr-Universit\"at Bochum,
D-44780 Bochum, Germany}
\email{rudolf.tange@rub.de}

\begin{document}

\begin{abstract}
 Let $G$ be a reductive linear algebraic group over an algebraically
closed field of characteristic $p > 0$.  A subgroup of $G$ is said
to be separable in $G$ if its global and infinitesimal
centralizers have the same dimension.  We study the interaction
between the notion of separability and Serre's concept of $G$-complete
reducibility for subgroups of $G$.  A separability hypothesis appears in
many general theorems concerning $G$-complete reducibility.  We
demonstrate that some of these results fail without this
hypothesis.  On the other hand, we prove that if $G$ is a connected
reductive group and $p$ is very good for $G$, then any subgroup of $G$
is separable; we deduce that under these hypotheses on $G$, a subgroup
$H$ of $G$ is $G$-completely reducible provided $\Lie G$ is semisimple
as an $H$-module.

 Recently, Guralnick has proved that if $H$ is a reductive subgroup of $G$ and
 $C$ is a conjugacy class of $G$, then $C\cap H$ is a finite union of
 $H$-conjugacy classes.  For generic $p$ --- when certain extra
 hypotheses hold, including separability --- this follows from a
 well-known tangent space argument due to Richardson, but in general,
 it rests on Lusztig's deep result that a connected reductive group
 has only finitely many unipotent conjugacy classes.  We show that the
 analogue of Guralnick's result is false if one considers conjugacy
 classes of $n$-tuples of elements from $H$ for $n > 1$.
\end{abstract}

\maketitle

\tableofcontents

\section{Introduction}
\label{sec:intro}
Let $G$ be a reductive algebraic group over an algebraically closed
field $k$.  A subgroup $H$ of $G$ is said to be $G$-completely
reducible ($G$-cr) if whenever $H$ is contained in a parabolic
subgroup $P$ of $G$, $H$ is contained in some Levi subgroup of $P$.
(If $G$ is non-connected, then we have to modify this definition
slightly; see Subsection \ref{subsec:noncon}.)  In the special case $G =
\GL_n(k)$, $H$ is $G$-completely reducible if and only if $H$ acts
completely reducibly on $k^n$.  Thus $G$-complete reducibility
generalizes the notion of semisimplicity from representation theory,
\cite{serre1.5, serre1, serre2}.

Much effort has been made to study $G$-completely reducible subgroups
of reductive groups.  There are applications to finite groups of Lie
type \cite{LMS}.  Ideas from the theory of $G$-complete reducibility
play an important part in the study by Liebeck and Seitz
\cite{liebeckseitz0} of the reductive subgroups of the exceptional
simple groups.

In characteristic zero a subgroup $H$ of $G$ is $G$-completely
reducible if and only if $H$ is reductive (cf.\ \cite[Lem.\ 2.6]{BMR}).  
In this case, it is
trivial to show that the class of $G$-completely reducible subgroups
of $G$ is closed under certain natural constructions: for instance, if
$H,N$ are $G$-completely reducible subgroups of $G$ such that $H$
normalizes $N$, then $H$ and $N$ are reductive, so $HN$ is reductive,
and thus $HN$ is $G$-completely reducible.  The analogous result need
not hold in positive characteristic \cite[Examples 5.1 and 5.3]{BMR2}.
Nonetheless, results concerning $G$-completely reducible subgroups 
which hold in
characteristic zero tend to hold ``generically'' in positive characteristic,
that is, if the
characteristic is large enough, or under certain extra natural
restrictions on the subgroups considered.  We give some other examples
below: see Theorem~\ref{thm:verygood}
and Theorem~\ref{thm:adjsc}.

The purpose of this paper is to explore the bounds of these generic
results, proving that they hold without the extra hypotheses or
finding counterexamples when the extra hypotheses are removed.
Characteristic $2$ is a particularly fertile ground for
counterexamples.  To prove positive results, we use geometric
techniques from our earlier works \cite{BMR} and \cite{BMR2}.
Especially important is the property of {\em separability} (see
\cite[Def.\ 3.27]{BMR}): the hypothesis that certain subgroups are
separable is required for several results involving $G$-complete
reducibility, e.g., \cite[Thm.\ 3.35, Thm.\ 3.46]{BMR}.

Let $\gg = \Lie G$ be the Lie algebra of $G$.
\begin{defn}
\label{def:sep}
 A subgroup $H$ of $G$ is said to be \emph{separable in $G$} if $\Lie
C_G(H) = \cc_\gg(H)$, that is, if the scheme-theoretic centralizer of $H$
in $G$ is smooth.
\end{defn}
Here $C_G(H)$ and $\cc_\gg(H)$
denote the (set-theoretic) centralizer of $H$ in $G$ and
in $\gg$ respectively.
Note that we always have
$\Lie C_G(H) \subseteq \cc_\gg(H)$.

The motivation for this terminology is as follows.  Given $n\in
\mathbb{N}$, we let $G$ act on $G^n$ by simultaneous conjugation:
\[
g\cdot (g_1,g_2,\ldots, g_n) = (gg_1g^{-1},gg_2g^{-1},\ldots,
gg_ng^{-1}).
\]
Suppose $H$ is the algebraic subgroup of $G$ generated by elements
$g_1, \ldots, g_n\in G$.  Then we say $H$ is \emph{(topologically)
finitely generated} by $g_1, \ldots, g_n$.  The orbit map $G \to
G\cdot(g_1,\ldots,g_n)$, $g\mapsto g\cdot (g_1,g_2,\ldots, g_n)$ is a
separable morphism of varieties if and only if
$H$ is separable in $G$ (see
\cite[Sec.~3.5]{BMR}).

There is an analogous notion of separability for subalgebras of $\gg$ ---
see Definition \ref{def:sepalg} --- which similarly is related to
the separability of orbit maps $G \to G\cdot(x_1,\ldots,x_n)$,
where $x_i \in \gg$ and $G$ acts diagonally on $\gg^n$ by simultaneous
adjoint action.

Our first main result concerning this notion of separability is as follows.

\begin{thm}
\label{thm:sep2}
Let $G$ be connected reductive and suppose that $\Char k$ is very good for
$G$.  Then any subgroup of $G$ is separable in $G$ and any subalgebra
of $\gg$ is separable in $\gg$.
\end{thm}

The proof (see Section \ref{sec:redpairs}) depends on Theorem \ref{thm:redpairs} below and the
existence of reductive pairs under the given characteristic
restrictions, where for a reductive subgroup $H$ of a reductive group $G$,
we say that $(G,H)$ is a \emph{reductive pair} if
the Lie algebra $\hh = \Lie H$ is an $H$-module direct
summand of $\gg$ (of course, this is automatically satisfied in
characteristic zero).
Special cases of Theorem \ref{thm:sep2}  are known due to work by
Liebeck--Seitz \cite[Thm.\ 3]{liebeckseitz0} and Lawther--Testerman
\cite[Thm.\ 2]{lawthertesterman}.

A prototype of the results we consider is the following fundamental
result of P.\ Slodowy \cite[Thm.~1]{slodowy}, which is obtained by applying a standard
tangent space argument due to R.W.\ Richardson
%\cite[Thm.~4.1]{rich2}.
(cf.\ \cite[Sec.~3]{rich2} and \cite[Lem.~3.1]{rich0}).
Theorem \ref{thm:finorbit} has many applications: see \cite{BMR},
\cite{slodowy}, \cite[Sec.~3]{vinberg}, for example.  Some of the
results that follow from Theorem \ref{thm:finorbit} are still valid
even when the hypotheses of the theorem are not met, but one often has
to work much harder to obtain them (cf.\ \cite[Sec.~1]{martin1}).

\begin{thm}
\label{thm:finorbit}
Let $H$ be a reductive subgroup of a reductive group $G$.
Let $n\in \mathbb{N}$, let $(h_1,\ldots,h_n)\in H^n$ and let $K$ be
the algebraic subgroup of $H$ generated by $h_1,\ldots,h_n$.  Suppose
that $(G,H)$ is a reductive pair and that $K$ is separable in $G$.
Then
\begin{itemize}
\item[(a)] for all $(g_1,\ldots,g_n)\in G\cdot (h_1,\ldots,h_n)\cap H^n$, the $H$-orbit map $H\ra H\cdot (g_1,\ldots,g_n)$ is
separable.  In particular, $K$ is a separable subgroup of $H$;
\item[(b)] the intersection $G\cdot (h_1,\ldots,h_n)\cap H^n$ is a
finite union of $H$-conjugacy classes;
\item[(c)] each $H$-conjugacy class in $G\cdot (h_1,\ldots,h_n)\cap H^n$ is
closed in $G\cdot(h_1,\ldots,h_n)\cap H^n$.
\end{itemize}
\end{thm}

The following consequence \cite[Thm.\ 3.35]{BMR} of Theorem
\ref{thm:finorbit} gives an application to $G$-complete reducibility.
The second assertion follows from Theorem \ref{thm:finorbit}(c) and
the geometric characterization of $G$-complete reducibility given in
\cite{BMR} (see Theorem \ref{thm:gcrcrit} below); the first does not
appear explicitly in \cite[Thm.\ 3.35]{BMR}, but it follows
immediately from the proof together with Theorem
\ref{thm:finorbit}(a).

\begin{thm}
\label{thm:redpairs}
Suppose that $(G,H)$ is a reductive pair. Let $K$ be a
subgroup of $H$ such that $K$ is a separable subgroup of $G$.  Then
$K$ is separable in $H$. Moreover, if $K$ is $G$-completely
reducible, then it is also $H$-completely reducible.
\end{thm}

Theorem \ref{thm:finorbit}(b) can be used to give a simple proof that
a connected reductive group $G$ has only finitely many conjugacy
classes of unipotent elements in generic characteristic: one takes an
embedding of $G$ in some $\GL(V)$ such that $(\GL(V),G)$ is a
reductive pair then applies Theorem \ref{thm:finorbit}(b) (taking $n=1$) to deduce
the result for $G$ from the result for $\GL(V)$.  In small positive
characteristic there need not exist any such reductive pair $(\GL(V),G)$ (see
Subsection \ref{subsec:sep}).  Nonetheless Lusztig proved without any
restriction on the characteristic, using some very deep mathematics,
that $G$ has only finitely many unipotent conjugacy classes
\cite{lusztig}.  Guralnick recently extended this result to
non-connected reductive groups \cite[Thm.\ 3.3]{guralnick}; he then
used it to prove that Theorem
\ref{thm:finorbit}(b) still holds for an arbitrary reductive subgroup $H$ of $G$ if $n=1$ \cite[Thm.\
1.2]{guralnick}, even though Richardson's proof no longer goes
through.

In view of this, it is natural to ask whether Theorem
\ref{thm:finorbit}(b) holds in more generality.  In Example
\ref{exmp:orbitcountereg}, we prove the following result which shows
that Guralnick's conjugacy result \cite[Thm.\ 1.2]{guralnick} does not
extend to conjugacy classes of $n$-tuples for $n>1$.

\begin{thm}
\label{thm:countereg}
 There exists a reductive group $G$, a reductive subgroup $M$ of $G$
 and a pair $(m_1,m_2)\in M^2$ such that $(G,M)$ is a reductive
 pair but $G\cdot (m_1,m_2)\cap M^2$ is {\bf not} a finite union
 of $M$-conjugacy classes.
\end{thm}

Thus Theorem \ref{thm:finorbit}(b) is false in general without the
separability hypothesis.  We give an example (Proposition
\ref{prop:GcrnotMcr}) showing that the second assertion of Theorem
\ref{thm:redpairs} also fails without the separability hypothesis.
This implies that Theorem
\ref{thm:finorbit}(c) fails as well without the separability
hypothesis (see Remark \ref{rem:thm1.3(c)fails}).  \smallskip

The paper is split into several sections, as we now outline.  In
Section \ref{sec:prelims} we introduce the notation and preliminary
results needed for the rest of the exposition.  In particular, we
recall the formalism of R-parabolic subgroups from \cite{BMR}, which
allows us to consider reductive groups which are not connected; this
is very important for many of our results.

In Section \ref{sec:redpairs} we discuss the question of separability
in connection with reductive pairs.  In Proposition
\ref{prop:C_G(S)redpair} we give a construction for certain reductive
pairs, where the reductive subgroup of $G$ is of the form $C_G(S)$; here $S$ is a reductive group acting suitably on $G$.

In Section~\ref{sec:sep1} we investigate the connection between
the $G$-complete reducibility of a subgroup $H$ and the semisimplicity of the adjoint
representation of $H$ on $\gg$.  The
following result is the basis of our discussion \cite[Thm.\
3.46]{BMR}.

\begin{thm}
\label{thm:adjcr}
Let $H$ be a separable subgroup of $G$.  If $\gg$ is semisimple as an
$H$-module, then $H$ is $G$-completely reducible.
\end{thm}

One way of removing the separability hypothesis from
Theorem~\ref{thm:adjcr} is to combine Theorems \ref{thm:sep2} and
\ref{thm:adjcr}. This immediately gives the following result.

\begin{thm}
\label{thm:verygood}
Let $G$ be connected reductive and suppose that $\Char k$ is very good
for $G$. Let $H$ be a subgroup of $G$ such that $\gg$ is a semisimple
$H$-module.  Then $H$ is $G$-completely reducible.
\end{thm}

The assumption in Theorem \ref{thm:verygood}
that $\Char k$  is very good for $G$ is rather
restrictive. In Section~\ref{sec:sep1} we discuss to what extent we
can remove the separability hypothesis from Theorem~\ref{thm:adjcr}
%when we only assume that $\Char k$  is good;
with a weaker assumption on $\Char k$;
see Theorem~\ref{thm:allregss}, Theorem~\ref{thm:adjsc} --- which is our second main result ---
and Corollary~\ref{cor:adjc}.

In Section \ref{sec:orbitthm} we extend two general results of Richardson
concerning orbits of reductive groups on affine varieties
\cite[Thm.\ A, Thm.\ C]{rich0}, see Theorem \ref{thm:richC}.
We apply these extensions in turn to
questions of $G$-complete reducibility; in particular,
we discuss some results which examine the
relationship between $G$-complete
reducibility and $H$-complete reducibility of a subgroup $K$ of $H$,
where $H$ is a subgroup of $G$ of the form $H = C_G(S)$,
with $S$ a reductive group acting suitably on $G$:
see Proposition \ref{prop:SnormalizesK}.
The case when the group $S$ considered above is a subgroup of $G$
is discussed in Section \ref{sec:cent+norm}.

In Section \ref{sec:ex}, we provide an important collection of closely
related constructions.  They give counterexamples to several of our
results, including Theorems \ref{thm:finorbit} and
\ref{thm:redpairs}, under weakened hypotheses.
In particular, here we also prove Theorem \ref{thm:countereg}.

\section{Preliminaries}
\label{sec:prelims}

\subsection{Notation}
Throughout, we work over an algebraically closed field $k$;
we let $k^*$ denote the multiplicative
group of $k$.
All algebraic groups are assumed to be linear.  By a subgroup of an algebraic group we mean a closed subgroup
and by a homomorphism of algebraic groups we mean a homomorphism
of abstract groups that is also a morphism of
algebraic varieties.  Let $H$ be a linear algebraic group.  We denote
by $\overline{\langle S\rangle}$ the algebraic subgroup of $H$
generated by a subset $S$.
We let $DH$ denote the derived group $[H,H]$, $Z(H)$ the centre of $H$,
and $H^0$ the connected component of $H$ that contains $1$.  If $S$ is
a subset of $H$, then $C_H(S)$ is the centralizer of $S$ in $H$ and
$N_H(S)$ is the normalizer of $S$ in $H$.  In general we use
an upper-case roman letter $G,H,K$, etc., to denote an algebraic group
and the corresponding lower-case gothic letter $\gg,\hh,\kk$, etc., to
denote its Lie algebra. If $\hh$ is a Lie algebra and $S$ is a
subset of $\hh$, then $\cc_\hh(S)$ is the centralizer of $S$ in
$\hh$.  We denote the centre of $\hh$ by $\zz(\hh)$.

Let $\Ad: H \to \GL(\hh)$ denote the adjoint representation; then we
let $H_{\rm ad}$ denote the image of $H$ under this map and $\hh_{\rm
ad}$ denote $\Lie H_{\rm ad}$.  Note that $(H_{\rm ad})^0$ is the adjoint form of $D(H^0)$ \cite[V.24.1]{borel}.

For the set of cocharacters (one-parameter subgroups) of $H$ we write
$Y(H)$; the elements of $Y(H)$ are the homomorphisms from $k^*$ to
$H$.

The \emph{unipotent radical} of $H$ is denoted $R_u(H)$; it is the
maximal connected normal unipotent subgroup of $H$.  The algebraic
group $H$ is called \emph{reductive} if $R_u(H) = \{1\}$; note that we
do not insist that a reductive group is connected.  In particular, $H$
is reductive if it is simple as an algebraic group ($H$ is said to be
\emph{simple} if $H$ is connected and all proper normal subgroups of
$H$ are finite).  If $N$ is a normal subgroup of $H$, then $H$ is reductive if and only if $N$ and $H/N$ are. The algebraic group $H$ is called \emph{linearly
reductive} if all rational representations of $H$ are semisimple.

If $H$ acts on the affine variety $X$, then we denote by $X^H$ the
fixed point subvariety of $X$: that is, $X^H = \{x \in X \mid h\cdot x
= x \ \forall h \in H\}$.  If $S$ is a subset of $X$, then we denote the pointwise
stabilizer of $S$ in $H$ by $C_H(S)$; we write $C_H(x)$ instead of
$C_H(\{x\})$ for $x$ in $X$.
If $X = K$ is an algebraic group and $H$ acts on $K$ by
automorphisms, then we write $C_K(H)$ instead of $K^H$.
Then we also have an
induced linear action of $H$ on $\kk = \Lie K$; we write $\cc_\kk(H)$ instead
of $\kk^H$.

Throughout the paper $G$ denotes a reductive algebraic group, possibly
non-connected, with Lie algebra $\gg$.
A subgroup of $G$ normalized
by some maximal torus of $G$ is called a \emph{regular} subgroup
of $G$ (connected reductive regular subgroups of connected reductive
groups are often also referred to as \emph{subsystem subgroups}, e.g.,
see \cite{liebeckseitz}).

Fix a maximal torus $T$ of $G$.
We write $X(T)$ for the character group of $T$.
Let $\Psi = \Psi(G,T) \subseteq X(T)$
denote the set of roots of $G$ with respect to
$T$. We write $\frakt = \Lie T$ for the Lie algebra of $T$.
If $\alpha\in \Psi$, then $U_\alpha$ denotes the
corresponding root subgroup of $G$ and $\uu_\alpha$ denotes the root space
${\rm Lie}\,U_\alpha$ of $\gg$.
Thus the root space decomposition of $\gg$ is
given by
\[
\gg = \frakt \oplus \bigoplus_{\alpha \in \Psi} \uu_\alpha.
\]
We denote by $G_\gamma$ the simple rank $1$
subgroup $\langle U_\gamma\cup U_{-\gamma}\rangle$ of $G$ and by $\gg_\gamma$
the Lie algebra of $G_\gamma$. Fix a Borel subgroup $B$ of $G$
containing $T$ and let $\Sigma = \Sigma(G, T)$ be the set of simple
roots of $\Psi$ defined by $B$. Then $\Psi^+ = \Psi(B,T)$ is the set of
positive roots of $G$.
For $\beta \in \Psi^+$ write $\beta =
\sum_{\alpha \in \Sigma} c_{\alpha\beta} \alpha$ with $c_{\alpha\beta}
\in \mathbb N_0$.  A prime $p$ is said to be \emph{good} for $G$ if it
does not divide any non-zero $c_{\alpha\beta}$, and
\emph{bad} otherwise.  A prime $p$ is good for $G$ if and only if it
is good for every simple factor of $G^0$ \cite{SS}; the bad primes for
the simple groups are $2$ for all groups except type $A_n$, $3$ for
the exceptional groups and $5$ for type $E_8$.  A prime $p$ is said to
be \emph{very good} for $G$, if $p$ is good for $G$ and $p$ does not
divide $n+1$ for any simple component of $G$ of type $A_n$. If $G$ is
simple and $\Char k$ is very good for $G$, then the Lie algebra $\gg$ is
simple \cite{St3}.

\begin{rem}
Separability of subgroups of $G$ and of subalgebras of $\gg$
(see Definition \ref{def:sepalg}) is automatic
in characteristic zero (cf.\ \cite[Thm.~13.4]{Hum}). Likewise,
the notion of $G$-complete reducibility is not
interesting in characteristic zero,
as a subgroup of $G$ is $G$-completely reducible if and only if
it is reductive (cf.\ \cite[Lem.\ 2.6]{BMR}); most of our results
and proofs become trivial in characteristic zero.
In the remainder of the paper
$p$ denotes the characteristic $\Char k$ of $k$ in case $\Char k>0$.
\end{rem}

Let $\gamma\in \Psi$.  We denote by $\gamma^\vee\in Y(G)$ the
corresponding coroot.  Then $\gamma^\vee$ is a homomorphism from $k^*$
to $G_\gamma$.  If $\alpha,\beta\in \Sigma$, then we have
$s_\alpha\cdot \beta= \beta- \langle \beta,\alpha^\vee\rangle \alpha$
\cite[Lem.~7.1.8]{spr2}, where $s_\alpha$ is the reflection
corresponding to $\alpha$ in the Weyl group of $G$.

In Section~\ref{sec:ex} we need the following well-known result
(which is implicit for instance in \cite[6.5]{Jantzen1}).
For convenience we include a proof.

\begin{lem}
\label{lem:sc}
Assume that $G$ is connected.  If the derived group $DG$ of $G$ is
simply connected, then the same holds for any Levi subgroup of any
parabolic subgroup of $G$.
\end{lem}

\begin{proof}
 We use the following characterization of simply connectedness of $DG$:
Let $T$ be a maximal torus of $G$ and let $\alpha_1,\ldots,\alpha_n\in X(T)$
be the choice of simple roots corresponding to some Borel subgroup $B$ containing $T$. Then $DG$ is simply connected if and only if
there exist characters $\chi_1,\ldots,\chi_n\in X(T)$ such that
$\inprod{\chi_i,\alpha_j\check}=\delta_{ij}$ for all $i,j\in\{1,\ldots,n\}$.
For $G$ semisimple this is clear. The general case follows from the fact
that the integer $\inprod{\chi,\alpha\check}$ only depends on the
restriction of $\chi$ to $T\cap DG$ and the fact that any character
of $T\cap DG$ can be lifted to a character of $T$
(\cite[Prop.~III.8.2(c)]{borel}).

 Let $P$ be a parabolic subgroup of $G$ and let $L$ be a Levi subgroup of $P$.  We may assume that $P$ contains $B$ and $L$ contains $T$.  The result now follows immediately from the following
well-known description of $L$
\cite[Prop.~IV.14.18]{borel}: the simple roots of $L$ can
be chosen from the set $\{\alpha_1,\ldots,\alpha_n\}$.
\end{proof}

The next result allows us in positive characteristic
to replace an algebraic group $S$ acting on $G$ by automorphisms
with a finite subgroup of $S$. It is a slight strengthening of \cite[Thm.~7]{BGM},
since it asserts the existence of a {\it countable} locally finite dense subgroup.
We use this lemma in the proof of
Proposition~\ref{prop:C_G(S)redpair}
and also in several places in Section~\ref{sec:orbitthm}.

\begin{lem}
\label{lem:finitedense}
Assume that $\Char k  > 0$.
Let $H$ be a linear algebraic group over $k$.
Then there exists an ascending sequence
$H_1\subseteq H_2\subseteq\cdots$ of finite
subgroups of $H$ whose union is dense in $H$.
\end{lem}

\begin{proof}
We proceed by induction on $\dim H$.
If $H$ is reductive, then the result follows
from \cite[Sec.~3]{martin1}.
Otherwise $Z:=Z(R_u(H))^0$ is a connected
unipotent normal subgroup of $H$ of dimension $\ge 1$.
By \cite[III.10.6(2)]{borel}, $Z$ contains a subgroup
isomorphic to the additive group $\mathbb{G}_a$.
Let $C$ be the subgroup of $Z$
generated by the subgroups of $Z$ that are
isomorphic to $\mathbb{G}_a$. Then $C$ is the additive group
of a non-zero finite-dimensional vector space over $k$,
by \cite[Thm.~5.4]{Hoch}.
%A connected abelian unipotent group need not be a vector space.
%See \cite{Hoch} Chap. 6 Sect. 5 note 2 p92.
Furthermore, $C$ is normal in $G$.
Clearly, $C$ has an ascending sequence $C_1\subseteq C_2\cdots$
of finite subgroups whose union is dense in $C$. By the
induction hypothesis, $M:=H/C$ also has an ascending sequence
$M_1\subseteq M_2\subseteq\cdots$ of finite subgroups whose
union is dense in $M$. Let $\pi:H\to M$ be the canonical
projection. For each $i\ge1$ let $H_i$ be a finitely generated
subgroup of $H$ such that $\pi(H_i)=M_i$. Without loss of
generality we may assume that the $H_i$ form an ascending
sequence of subgroups and that $H_i$ contains $C_i$. Since
$H_i$ is finitely generated and $H_i\cap C$ is of finite index
in $H_i$, we have that $H_i\cap C$ is finitely generated,
by \cite[Thm.~11.54]{Rot}.
%The cited theorem say that a subgroup of finite index of a free group
%of finite rank is free of finite rank. The condition "of finite index"
%is necessary.
%This clearly implies that any subgroup of finite index of a finitely
%generated group is finitely generated.
%The Schreier-Nielsen Theorem says that a subgroup of a free group is free
%which is not exactly what we need.
Since $C$ is a vector space, this means that $H_i\cap C$ is
finite. But then $H_i$ is finite. Now let $H'$ be the closure of
the union of the $H_i$. Then $H'$ is a closed subgroup of $H$
containing $C$. Its image $\pi(H')$ is a closed subgroup of $M$
containing the $M_i$ and is therefore equal to $M$. Consequently, $H'=H$.
\end{proof}

\subsection{$G$-Complete Reducibility}
\label{subsec:noncon}
In \cite[Sec.~6]{BMR}, Serre's original notion of $G$-complete
reducibility is extended to include the case when $G$ is reductive but
not necessarily connected (so that $G^0$ is a connected reductive
group).  The crucial ingredient of this extension is the introduction
of so-called \emph{Richardson-parabolic subgroups} (\emph{R-parabolic
subgroups}) of $G$.  We briefly recall the main definitions here; for
more details on this formalism, see \cite[Sec.~6]{BMR}.
%and \cite[Sec.~2]{BMR2}.

For a cocharacter $\lambda
\in Y(G)$, the \emph{R-parabolic subgroup} corresponding to $\lambda$
is defined by $P_\lambda := \{ g\in G \mid \underset{a\to 0}{\lim}\,
\lambda(a)g\lambda(a)\inverse \textrm{ exists}\}$.  Then $P_\lambda$
admits a Levi decomposition $P_\lambda = R_u(P_\lambda) \rtimes
L_\lambda$, where $L_\lambda = \{ g\in G \mid \underset{a\to
0}{\lim}\, \lambda(a)g\lambda(a)\inverse = g \} = C_G(\lambda(k^*))$.
We call $L_\lambda$
an \emph{R-Levi subgroup} of $P_\lambda$.  For an R-parabolic subgroup
$P$ of $G$, the different R-Levi subgroups of $P$ correspond in this
way to different choices of $\lambda \in Y(G)$ such that $P =
P_\lambda$; moreover, the R-Levi subgroups of $P$ are all conjugate
under the action of $R_u(P)$.  An R-parabolic subgroup $P$ is a
parabolic subgroup in the sense that $G/P$ is a complete variety; the
converse is true when $G$ is connected, but not in general
(\cite[Rem.\ 5.3]{martin1}).  The map $c_\lambda :P_\lambda \to
L_\lambda$ given by $c_\lambda(g) = \underset{a\to 0}{\lim}\,
\lambda(a)g\lambda(a)\inverse$ is a surjective homomorphism of
algebraic groups with kernel $R_u(P_\lambda)$; it coincides with the
usual projection $P_\lambda\ra L_\lambda$.

\begin{rem}
For a subgroup $H$ of $G$, there is a natural inclusion $Y(H)
\subseteq Y(G)$.  If $\lambda \in Y(H)$, and $H$ is reductive, we can
therefore associate to $\lambda$ an R-parabolic subgroup of $H$ as
well as an R-parabolic subgroup of $G$.  To avoid confusion, we
reserve the notation $P_\lambda$ for R-parabolic subgroups of $G$, and
distinguish the R-parabolic subgroups of $H$ by writing $P_\lambda(H)$
for $\lambda \in Y(H)$.  The notation $L_\lambda(H)$ has the obvious
meaning.  Note that $P_\lambda(H) = P_\lambda \cap H$,
$L_\lambda(H) = L_\lambda \cap H$ and
$R_u(P_\lambda(H))= R_u(P_\lambda)\cap H$ for $\lambda \in Y(H)$.
\end{rem}

\begin{rem}
\label{rem:P0G0}
If $\lambda\in Y(G)$ and $P_\lambda^0=G^0$,
then $P_\lambda$ is an R-Levi subgroup of itself:
for $R_u(P_\lambda)$, being connected, is contained in $G^0$, so
$R_u(P_\lambda)= R_u(P_\lambda)\cap G^0= R_u(P_\lambda(G^0))= R_u(P_\lambda\cap G^0)= R_u(G^0)= \{1\}$.
\end{rem}

\begin{defn}
\label{defn:gcr}
Suppose $H$ is a subgroup of $G$.  We say $H$ is \emph{$G$-completely
reducible} ($G$-cr for short) if whenever $H$ is contained in an
R-parabolic subgroup $P$ of $G$, then there exists an R-Levi subgroup
$L$ of $P$ with $H \subseteq L$.
\end{defn}

\begin{rem}
If $H$ is a $G$-completely reducible subgroup of $G$,
then $H$ is reductive (cf.\ \cite[Sec.~2.5 and Thm.~3.1]{BMR}).
\end{rem}

\begin{comment}
We give a useful characterisation of $G$-complete reducibility in
terms of the map $c_\lambda$ \cite[Lem.~2.17 and Thm.~3.1]{BMR}.

\begin{lem}
\label{lem:gcrclambda}
 Let $H$ be a subgroup of $G$.  Then $H$ is $G$-completely reducible
 if and only if for every $\lambda\in Y(G)$ such that $H\subseteq
 P_\lambda$, $c_\lambda(H)$ is $R_u(P_\lambda)$-conjugate to $H$.
\end{lem}
\end{comment}

Since all parabolic subgroups (respectively\ all Levi subgroups of
parabolic subgroups) of a connected reductive group are R-parabolic
subgroups (respectively\ R-Levi subgroups of R-parabolic subgroups),
Definition \ref{defn:gcr} coincides with Serre's original definition
for connected groups \cite{serre2}.
Sometimes we come across subgroups of $G$ which
are not contained in any proper R-parabolic subgroup of $G$; these
subgroups are trivially $G$-completely reducible.  Following Serre
again, we call these subgroups \emph{$G$-irreducible} ($G$-ir).
%Similarly, a subgroup of $G$ which is not contained in any R-Levi
%subgroup of any proper R-parabolic subgroup is called
%\emph{$G$-indecomposable} ($G$-ind).

\begin{rem}
\begin{comment}
 If $H$ is a $G$-completely reducible subgroup of $G$ and $P$ is
minimal among R-parabolic subgroups of $G$ containing $H$, then $H$ is
$L$-irreducible, where $L$ is an R-Levi subgroup of $P$ containing
$H$.  This observation often allows us to restrict attention to
$G$-irreducible and $G$-indecomposable subgroups.
\end{comment}
Since R-Levi
subgroups of R-parabolic subgroups play an important r\^ole in many of
our proofs, for brevity we sometimes abuse language and refer to an
\emph{R-Levi subgroup of $G$}; by this we mean an R-Levi subgroup of
some R-parabolic subgroup of $G$.  Similarly, when $G$ is connected,
we may refer to \emph{a Levi subgroup of $G$}; this means a Levi
subgroup of some parabolic subgroup of $G$.
\end{rem}

A key result is the following \cite[Cor.\ 3.7]{BMR}, which gives a
geometric criterion for $G$-complete reducibility.

\begin{thm}
\label{thm:gcrcrit}
 Let $g_1,\ldots, g_n\in G$ and let $H= \overline{\langle \{g_1,\ldots,
 g_n\}\rangle}$.  Then $H$ is $G$-completely reducible if and only if
 the conjugacy class $G\cdot (g_1,\ldots, g_n)$ is closed in $G^n$.
\end{thm}

%Proposition \ref{prop:LsepGnonsep} shows that the converse of the %first assertion of Theorem \ref{thm:redpairs}
%is false in general.

We frequently require results from \cite[Sec.~6.3]{BMR} for non-connected $G$,
though
we usually simply cite the relevant result in \cite{BMR} for connected $G$.

\subsection{Separability}
\label{subsec:sep}
In Section~\ref{sec:redpairs} we require the following analogue of
Definition~\ref{def:sep} for subalgebras of $\gg$.
Recall that
$C_G(\hh) = \{g \in G \mid \Ad g(x) = x \text{ for all } x \in \hh \}$.
%The notation
%$C_G(\hh)$ was introduced in 2.1: take $H=G$, $X=\gg$ as the
%$G$-variety and $A=\hh$.

\begin{defn}
\label{def:sepalg}
A subalgebra $\hh$ of $\gg$ is said to be \emph{separable in $\gg$} if
$\Lie C_G(\hh) = \cc_\gg(\hh)$.
\end{defn}

The above definition has the same motivation as in the group
case. Given $n\in \mathbb{N}$, we let $G$ act on $\gg^n$ by diagonal
adjoint action.  Suppose $\hh$ is the subalgebra of $\gg$ generated by
elements $x_1, \ldots, x_n\in \gg$. Then the orbit map $G \to
G\cdot(x_1,\ldots,x_n)$, $g\mapsto g\cdot (x_1, \ldots,x_n)$
is a separable morphism of varieties if and only if $\hh$ is separable in
$\gg$ (see \cite[II.6.7]{borel}).

As in the group case (\cite[Thm.~1]{slodowy}), it is straightforward
to generalize Richardson's tangent space arguments
%\cite[Thm.\ 4.1]{rich2}
\cite[Sec.~3]{rich2} and \cite[Lem.~3.1]{rich0}
to the action of $G$ on $\gg^n$. We
then obtain the following analogue of Theorem~\ref{thm:finorbit}.

%{\tt R: Do we have a reference for this like in the group case? Maybe
%this has been stated before. We could refer to this theorem with a
%phrase like: "We obtain the analogue of ...", rather than stating it
%explicitly.}

\begin{thm}
\label{thm:finorbitalg}
Let $H$ be a reductive subgroup of $G$.
Let $n\in \mathbb{N}$, let $(x_1,\ldots,x_n)\in \hh^n$ and let $\kk$
be the subalgebra of $\hh$ generated by $x_1,\ldots,x_n$.  Suppose
that $(G,H)$ is a reductive pair and that $\kk$ is separable in $\gg$.
Then
\begin{itemize}
\item[(a)] all $H$-orbit maps in $G\cdot (x_1,\ldots,x_n)\cap \hh^n$
are separable. In particular, $\kk$ is separable in $\hh$;
\item[(b)] the intersection $G\cdot (x_1,\ldots,x_n)\cap \hh^n$ is a
finite union of $H$-conjugacy classes;
\item[(c)] each $H$-conjugacy class in $G\cdot (x_1,\ldots,x_n)\cap \hh^n$ is
closed in $G\cdot(x_1,\ldots,x_n)\cap \hh^n$.
\end{itemize}
\end{thm}

We immediately obtain the following analogue of the first
assertion of Theorem~\ref{thm:redpairs}.

\begin{cor}
\label{cor:redpairs}
Suppose that $(G,H)$ is a reductive pair.  Let $\kk$ be a subalgebra
of $\hh$ such that $\kk$ is separable in $\gg$.  Then $\kk$ is
separable in $\hh$.
\end{cor}

Because every subgroup of $\GL(V)$ is separable in $\GL(V)$ and every
subalgebra of $\gl(V)$ is separable in $\gl(V)$ (e.g., see \cite[Ex.\
3.28]{BMR}), Theorem~\ref{thm:redpairs} and
Corollary~\ref{cor:redpairs} imply the following.

\begin{cor}
\label{cor:glv:g}
If $(\GL(V),G)$ is a reductive pair, then every subgroup of $G$ is
separable in $G$ and every subalgebra of $\gg$ is separable in $\gg$.
\end{cor}

We deduce from this that not every $G$ can be embedded in some
$\GL(V)$ in such a way that $(\GL(V),G)$ is a reductive pair: this
applies, for instance, to $G=\SL_2(k)$ when $p=2$, because $H=G$ is
not a separable subgroup of $G$.  However, if $G$ is of a given Dynkin
type, then generically --- that is, for almost all values of $p$ ---
the conclusion of Corollary \ref{cor:glv:g} holds; for example, if $G$
is an exceptional simple group of adjoint type and $p$ is good for
$G$, then $(\GL(\gg),G)$ is a reductive pair (cf.\ Example
\ref{exmp:Serre}).

\medskip

The final result of this section shows that a non-separable $G$-cr
subgroup $K$ of $G$ is, up to isogeny, a separable subgroup of a
regular subgroup of $G$.
Given a reductive group $M$, we let $\pi_M : M\ra M_{\rm ad}$ denote the natural morphism.

\begin{prop}
\label{prop:sepovergp}
Let $K$ be a $G$-completely reducible subgroup of $G$.  Then there
exists a reductive subgroup $M$ of $G$ containing a maximal torus of $G$
%regular reductive subgroup $M$ of $G$
such that $K \subseteq M$, $K$ is $M$-irreducible and $\pi_M(K)$ is
separable in $M_{\rm ad}$.
\end{prop}

\begin{proof}
Since
%a regular reductive subgroup of an R-Levi subgroup of $G$ is a
%regular reductive subgroup of $G$ and
$K$ is $G$-cr, we may assume by \cite[Cor.~3.5]{BMR}
that $K$ is $G$-ir after replacing $G$ by an R-Levi subgroup of $G$
that is minimal with respect to containing $K$.  By
\cite[Cor.~2.7(i)]{BMR}, $K$ is $M$-ir in any reductive subgroup $M$
of $G$ containing $K$.  If $\pi_G(K)$ is separable in $G_{\rm ad}$,
then we can take $M=G$, so suppose not.
By \cite[Prop.\ 3.39]{BMR}, there exists a %regular
reductive subgroup $M'$ of $G_{\rm ad}$ containing a maximal torus of $G_{\rm
ad}$ such that $\pi_G(K)\subseteq M'$ and $M'$ is not separable in
$G_{\rm ad}$.  As $(G_{\rm ad})^0$ is of adjoint type, its
Lie algebra has trivial centre,
so any overgroup of $(G_{\rm ad})^0$ is separable in $G_{\rm ad}$.
This forces $M'$ to be of strictly smaller dimension than $G_{\rm
ad}$.  Let $M= \pi_G^{-1}(M')$, a subgroup of $G$ which is of strictly
smaller dimension than $G$ and contains a maximal torus of $G$.  Since
$M$ is an overgroup of the $G$-ir group $K$, $M$ is reductive.  The
result now follows by induction on $\dim G$.
\end{proof}

Proposition \ref{prop:sepovergp} is false if we do not assume that $K$
is $G$-completely reducible: see Example \ref{exmp:nonGcr} below.

\section{Reductive Pairs and Separability}
\label{sec:redpairs}
The notion of separability is central to many of the results in this
paper.  Theorems \ref{thm:finorbit} and \ref{thm:redpairs} both
illustrate the importance of reductive pairs in this context.  In this
section we elaborate on this theme.  For examples and constructions of
reductive pairs, we refer to \cite[Sec.~I.3]{slodowy} and \cite[Sec.~3.5]{BMR}.

%The following lemma is used in the sequel.
Recall that an isogeny
is an epimorphism with finite kernel and that it
is called \emph{separable}
if its differential is an isomorphism.

\begin{lem}
\label{lem:sep}
Let $\varphi:G\to G'$ be a separable isogeny of reductive
groups, let $H$ be a subgroup of $G$ and let $\kk$ be a subalgebra of
$\gg$.  Then $H$ is separable in $G$ if and only if $\varphi(H)$ is
separable in $G'$ and $\kk$ is separable in $\gg$ if and only if
$d\varphi(\kk)$ is separable in $\gg'$.
\end{lem}
\begin{proof}
Since the differential $d\varphi$ of $\varphi$ is an isomorphism, it is
clear that
$d\varphi(\cc_\gg(\kk))=\cc_{\gg'}(d\varphi(\kk))$. Furthermore, we
have
 \begin{equation}
 \label{eqn:ad}
\Ad(\varphi(g))\circ d\varphi
=d\varphi\circ\Ad(g)\text{\quad for all } g\in G,
 \end{equation}
from which it follows that
$d\varphi(\cc_\gg(H))=\cc_{\gg'}(\varphi(H))$.

So to prove the first statement it suffices to show that
$\varphi(C_G(H))^0=C_{G'}(\varphi(H))^0$.  For this, in turn, it
suffices to show that $\varphi^{-1}(C_{G'}(\varphi(H)))^0\subseteq
C_G(H)$.  If $y\in H$, then the image of
$\varphi^{-1}(C_{G'}(\varphi(H)))^0$ under the morphism $x\mapsto
xyx^{-1}y^{-1}$ is an irreducible subset of the finite set $\ker \varphi$ which
contains $1$, so it must equal $\{1\}$.

To prove the second statement, it suffices to show that
$\varphi(C_G(\kk))=C_{G'}(d\varphi(\kk))$. This follows easily from
Eqn.\ \eqref{eqn:ad}.
\end{proof}

\begin{lem}
\label{lem:sepisog}
Let $G$ be connected and suppose that $p$ is very good for
$G$.  Then there exists a separable isogeny $S\times H\to G$, where
$S$ is a torus and $H$ is a product of simply connected simple groups.
\end{lem}

\begin{proof}
 Let $G_1,\ldots, G_r$ be the simple factors of $DG$ and let
 $\widetilde{G}_i$ be the simply connected cover of $G_i$ for each
 $i$.  Then $\Lie \widetilde{G}_i$ is simple for each $i$, by our
 hypothesis on $p$.  Set $S=Z(G)^0$ and $H= \widetilde{G}_1\times
 \cdots \times \widetilde{G}_r$.  It is easily checked that the
 multiplication map $S\times H\ra G$ is a separable isogeny.
\begin{comment}
By \cite{St2} there exists a simply connected semi-simple connected
group $H$ and a central isogeny $H\to DG$. Let $(X,\Phi)$ be the root
datum of $DG$ and let $X_r$ resp.\ $X_w$ be the root lattice resp.\
the weight lattice of $\Phi$. Then the kernel of the differential of
this isogeny is isomorphic to $\Hom(X_w/X,k)$. Now $X_w/X$ is a
quotient of $X_w/X_r$ and an easy inspection of \cite[Planches I -
IX]{bou} (one may assume for this that $\Phi$ is irreducible) shows
that a very good prime does not divide the order of $X_w/X_r$.  So the
differential of the isogeny has trivial kernel and the isogeny is
separable.

Now $H$ is the direct product of almost simple, simply connected
groups (see \cite{St2}) and the Lie algebras of these groups are
simple (see \cite{St3}, from our assumptions on the characteristic of
$k$ it follows that a basis of the root system of such a group gives
rise to a basis of the dual of the Lie algebra of a maximal torus), so
$\Lie DG$ is the direct sum of simple Lie algebras.  In particular,
$\Lie DG$ is perfect and centreless.  Since we always have $[\gg,\gg]
\subseteq \Lie DG$ (\cite[3.17]{borel}), it follows that $\Lie DG =
[\gg,\gg]$.  Furthermore we have $\zz(\gg) \cap [\gg,\gg] = \{0\}$,
since $[\gg,\gg]$ is centreless. The multiplication $Z(G)^0\times DG
\to G$ is a central isogeny and this isogeny is separable, since $\Lie
Z(G) \cap \Lie DG = \{0\}$.  Combining this with the isogeny $H\to DG$
we obtain the assertion of the lemma.
\end{comment}
\end{proof}

We are now in a position to prove Theorem \ref{thm:sep2}.

%\begin{thm}
%\label{thm:sep2}
%Let $G$ be connected reductive and suppose that $p$ is very good for
%$G$.  Then any subgroup of $G$ is separable in $G$ and any subalgebra
%of $\gg$ is separable in $\gg$.
%\end{thm}

\begin{proof}[Proof of Theorem \ref{thm:sep2}]
By Lemmas~\ref{lem:sep} and \ref{lem:sepisog},
we may assume that $G =
S\times H_1\times \cdots \times H_r$, where $S$
is a torus and each $H_i$ is a simply connected
simple group. Put $H_0=S$ and let $\pi_i:G\to H_i$ be the projection.
For every subgroup $K$ of $G$ we have
$C_G(K)=\prod_{i=0}^rC_{H_i}(\pi_i(K))$ and for every subalgebra $\mm$ of $\gg$ we have
$\cc_\gg(\mm)=\bigoplus_{i=0}^r\cc_{\hh_i}(d\pi_i(\mm))$. Since $S$ is abelian, we may now
assume that $G$ is a simply connected simple group.

It now suffices to prove that there
exists a simple group $G'$, a separable isogeny
$\eta\colon G\ra G'$ and an embedding of $G'$ in
some $\GL(V')$ such that $(\GL(V'),G')$ is a reductive
pair: for then every subgroup of $G$ is separable in $G$
and every subalgebra of $\gg$ is separable in $\gg$, by
Corollary~\ref{cor:glv:g} and Lemma~\ref{lem:sep}.
For every simple group $K$ of the same Dynkin type as $G$,
the natural isogeny $G\to K$ is separable, since $p$ is very good for $G$.
So to complete the proof, it is enough to show
that for every Dynkin type, there exists a simple group
$K$ of this Dynkin type and an embedding of $K$ in some
$\GL(V)$ such that $(\GL(V),K)$ is a reductive pair.

If $K = \SO(V)$ or $K = \SP(V)$ and $p\neq 2$, then $(\GL(V),K)$ is a
reductive pair \cite[Lem.~5.1]{rich2}.  This deals with types $B_n$,
$C_n$ and $D_n$.  If $K = \SL(V)$, then, since $p$ is very good for
$G$, it follows that $(\GL(V),K)$ is a reductive pair: the scalar
matrices form a $K$-stable direct complement to $\Sl(V)$ in
$\gl(V)$. This deals with type $A_n$.  If $K$ is an adjoint simple
group of exceptional type and $p$ is good for $K$, then
$(\GL(\kk),K)$ is a reductive pair, thanks to \cite[\S 5]{rich2}.
This completes the proof.
\end{proof}

\begin{cor}
\label{cor:sep2}
Let $G$ be connected and suppose that $p$ is very good for
$G$.  Let $g_1, \ldots, g_n \in G$ and let $x_1, \ldots,x_m\in\gg$.
Then the orbit maps $G \to G \cdot (g_1, \ldots, g_n)$ and $G \to G
\cdot (x_1, \ldots, x_m)$
%\subseteq G^n$
are separable.
\end{cor}

\begin{proof}
 The orbit map $G \to G \cdot (g_1, \ldots, g_n)$ is separable if and
only if the algebraic subgroup of $G$ generated by $g_1, \ldots, g_n$
is separable in $G$ and the orbit map $G \to G\cdot(x_1,\ldots,x_m)$
is separable if and only if the subalgebra of $\gg$ generated by
$x_1,\ldots, x_m$ is separable in $\gg$ (see the comment
after Definition \ref{def:sep} and Subsection \ref{subsec:sep}),
so the corollary follows
immediately from Theorem \ref{thm:sep2}.
\end{proof}

\begin{rems}
\label{rems:redpair}
\begin{comment}
 \begin{itemize}
  \item[(i)] For $G$ simple of exceptional type and for simple
subgroups of $G$ and $p>7$, Theorem \ref{thm:sep2} is due to
case-by-case checks of Liebeck--Seitz \cite[Thm.\ 3]{liebeckseitz0}
and Lawther--Testerman \cite[Thm.\ 2]{lawthertesterman}.
  \item[(ii)] The restriction on $p$ in Theorem \ref{thm:sep2} is
necessary.  For instance, for $G = \SL(V)$ with $\dim V = p$, the
group $G$ is not separable in itself.  Also for $G$ simple of
exceptional type and $p$ a bad prime for $G$, the pair $(\GL(\gg), G)$
need no longer be a reductive pair.
%as there always exists a non-separable subgroup of $G$
For instance, in Proposition \ref{prop:LsepGnonsep} we provide an
example for $G$ of type $G_2$ and $p=2$ of a non-separable subgroup of
$G$ (cf.\ Corollary \ref{cor:glv:g}).
  \item[(iii)] The requirement that $G$ be connected in Theorem \ref{thm:sep2}
  is also necessary.  For instance, if $G=k^*\rtimes C_2$, where the
  non-trivial element $c$ of the cyclic group $C_2$ acts on $k^*$ by
  $c\cdot a = a^{-1}$, then $C_2$ is a non-separable subgroup of $G$.
  \item[(iv)] The case $n=1$ in Corollary \ref{cor:sep2} is a
  well-known fundamental result due to P.\ Slodowy, \cite[p.\ 38]{slodowy1}.
 \end{itemize}
\end{comment}
(i). Consider the class of reductive groups $G$ that have the property
that each subgroup of $G$ is separable in $G$ and each subalgebra of $\gg$
is separable in $\gg$. Lemma~\ref{lem:sep}, the proof of Theorem~\ref{thm:sep2}
and Corollary~\ref{cor:LGsep} below show that this class is closed under separable
isogenies (in both directions), direct products and centralizers of subgroups $S$ acting on $G$ by automorphisms
as in Proposition~\ref{prop:C_G(S)redpair}.
In particular our class contains the ``strongly standard'' reductive groups of
\cite[\S 2.4]{mcninch-test}.

(ii). For $G$ simple of exceptional type and for simple subgroups of
$G$ and $p>7$, Theorem \ref{thm:sep2} is due to case-by-case checks of
Liebeck--Seitz \cite[Thm.\ 3]{liebeckseitz0} and
Lawther--Testerman \cite[Thm.\ 2]{lawthertesterman}.

(iii). The restriction on $p$ in Theorem \ref{thm:sep2} is necessary.
For instance, for $G = \SL(V)$ with $\dim V = p$, the group $G$ is not
separable in itself.  Also for $G$ simple of exceptional type and $p$
a bad prime for $G$, the pair $(\GL(\gg), G)$ need no longer be a
reductive pair, so the proof breaks down:
%as there always exists a non-separable subgroup of $G$
for instance, in Proposition \ref{prop:LsepGnonsep} we provide an
example for $G$ of type $G_2$ and $p=2$ of a non-separable subgroup of
$G$ (cf.\ Corollary \ref{cor:glv:g}).

(iv). The requirement in Theorem \ref{thm:sep2} that $G$ be connected
is also necessary.  For instance, if $G=k^*\rtimes C_2$, where the
non-trivial element $c$ of the cyclic group $C_2$ acts on $k^*$ by
$c\cdot a = a^{-1}$, then $C_2$ is a
non-separable subgroup of $G$.

(v). The case $n=1$ in Corollary \ref{cor:sep2} is a well-known
fundamental result due to P.\ Slodowy, \cite[p38]{slodowy1}.

(vi). Serre has asked whether Theorem \ref{thm:sep2} holds for an arbitrary group subscheme $H$ of $G$;
Theorem \ref{thm:sep2} deals with the two special cases $H$ smooth and $H$ infinitesimal of height one.
\end{rems}

%Combined with Theorem~\ref{thm:redpairs},
%the following result shows that separability of a subgroup $H$ in %$G$
%implies separability of $H$ in any Levi subgroup of $G$
%containing $H$.
We finish the section with some further useful results on
separability and reductive pairs.

\begin{lem}
\label{lem:almostfaithful}
 Suppose $G$ is connected.
Let $S$ be an algebraic group
acting faithfully on $G$ by automorphisms.
Then the corresponding representation of $S$ on $\gg$ has finite kernel.
\end{lem}

\begin{proof}
 It is enough to prove that $S=\{1\}$ under the
extra hypotheses that $S$ is connected and $S$
acts trivially on $\gg$, so we assume this.
Let $s\in S$.  Let $B,B^-$ be any pair of
opposite Borel subgroups of $G$.  Since
$s\cdot B$ is also a Borel subgroup of $B$,
there exists $g\in G$ such that $s\cdot B= gBg^{-1}$.
Let $\bb= \Lie B$.  Then $\bb= s\cdot \bb= \Ad g(\bb)$,
so $g\in B$, by \cite[IV.14.1 Cor.~2]{borel}, so $s$
normalizes $B$.  Similarly $s$ normalizes $B^-$, so
$s$ normalizes the maximal torus $T=B\cap B^-$ of $G$.
Hence $s$ centralizes $T$, by the rigidity of tori
\cite[Prop.~III.8.10]{borel}. Since $B$ and $B^-$
were arbitrary, $s$ centralizes the set of semisimple elements of $G$,
which is dense in $G$ as $G$ is reductive (e.g., see \cite[Thm.~IV.12.3(2)]{borel}).
Thus $s$ centralizes $G$. But the
$S$-action on $G$ is faithful, so $S=\{1\}$ as required.
\end{proof}

%{\tt R: Below we could omit the assumption that $S$ is reductive.
%See after the proof in the tex file. I think it isn't worth it.}

\begin{prop}
\label{prop:C_G(S)redpair}
Let $S$ be
%a reductive
an algebraic group acting on $G$ by automorphisms.
Suppose that $S$ acts semisimply on $\gg$ and $\Lie C_G(S)=\cc_\gg(S)$.
Then
 \begin{itemize}
  \item[(a)] $C_G(S)$ is $G$-completely reducible
  and $(G,C_G(S))$ is a reductive pair.
  \item[(b)] If $S$ is a subgroup of $G$, then $N_G(S)$ is $G$-completely reducible and
  if further $S^0$ is central in $S$, then $(G,N_G(S))$ is a
  reductive pair.
 \end{itemize}
\end{prop}

\begin{proof}
 (a). Let $H$ be the union of the components of $G$ that meet $C_G(S)$.  Then $H$ is a finite-index subgroup of $G$, so $C_G(S)$ is $G$-cr if and only if it is $H$-cr, by \cite[Prop.~2.12]{BMR2}.  Hence we can assume that $G=H$.  Replacing $S$ by $S/C_S(G)$, we can also assume that $S$ acts faithfully on $G$.  Since $C_G(S)$ meets every component of $G$, it follows that $S$ acts faithfully on $G^0$.  The completely reducible representation $S\ra \GL(\gg)$ has finite kernel by Lemma \ref{lem:almostfaithful}, so $S$ is reductive.

If $\Char k=0$, then, as $S$ is reductive, it is linearly reductive.
So in this case all we have to show is that $C_G(S)$ is reductive.
Since $C_G(S)^0=C_{G^0}(S)^0$, this follows immediately from
\cite[Prop.~10.1.5]{rich0}.

Now assume that $\Char k = p > 0$.
\begin{comment}
Since $S$ is reductive, $S$ has an
$\overline{\mathbb{F}_p}$-structure, thanks to \cite[Prop.~3.2]{martin1}.
As pointed out in \cite[Notation~3.3]{martin1}, this implies that $S$ has an
increasing chain of finite subgroups $S_1\subseteq S_2\subseteq\cdots$
such that $\bigcup_{i\ge1}S_i$ is Zariski dense in $S$.  From this it
follows that for some $i\ge1$, $\gg$ is a semisimple $S_i$-module (cf.\ \cite[Lem.~2.10]{BMR}),
and we have $C_G(S_i) = C_G(S)$ and
$\cc_\gg(S_i) = \cc_\gg(S)$.
\end{comment}
By Lemma~\ref{lem:finitedense}, $S$ admits an
ascending chain of finite subgroups $S_1\subseteq S_2\subseteq\cdots$
such that
%$\bigcup_{i\ge1}S_i$
their union is Zariski dense in $S$. From this it
follows that for some $i\ge1$, $\gg$ is a semisimple $S_i$-module (that is, the image of $S_i$ in $\GL(\gg)$ is $\GL(\gg)$-cr --- cf.\ \cite[Lem.\ 2.10]{BMR}) and we have $C_G(S_i) = C_G(S)$ and $\cc_\gg(S_i) = \cc_\gg(S)$.
So, after replacing $S$ by $S_i$, we may assume that $S$ is finite.
Put $G' = G \rtimes  S$. Clearly, $G'$ is reductive.
Furthermore, $\gg' = \Lie G' = \gg$ is a semisimple $S$-module,
and $\Lie C_{G'}(S) = \cc_{\gg'}(S)$ by construction.
Therefore, $S$ is $G'$-cr, by Theorem \ref{thm:adjcr}.
%see also \cite[Sec.~6.3]{BMR} last sentence.
Since $S$ is $G'$-cr, $C_{G'}(S)$ is $G'$-cr,
by \cite[Cor.\ 3.17]{BMR}.
So the normal subgroup $C_G(S)$ of $C_{G'}(S)$ is also $G'$-cr,
\cite[Thm.\ 3.10]{BMR}.
Now $C_G(S)$ is $G$-cr, by \cite[Prop.\ 2.12]{BMR2}.

\begin{comment}
%R:this is incorrect as Gerhard pointed out
%$C_{\GL(\gg)}(S)^0$ could be much bigger than the image of
%$C_G(S)^0$
We have a completely reducible representation from $S$ into
$\GL(\gg)$, so $C_{\GL(\gg)}(S)$ is $\GL(\gg)$-cr. Let $\pi_G\colon
G\ra \GL(\gg)$ be the canonical projection.  Then
$\pi_G(C_G(S))\subseteq C_{\GL(\gg)}(S)$, so $C_G(S)\subseteq
\pi_G^{-1}(C_{\GL(\gg)}(S))$.  Hence $C_G(S)^0$ is an extension of the
reductive group $C_{\GL(\gg)}(S)^0$ by a subgroup of the linearly
reductive group $Z(G^0)$.  This implies that $C_G(S)^0$ is reductive.
\end{comment}

Since $\gg$ is a semisimple $S$-module, it has a direct sum
decomposition into $S$-isotypic summands.
By hypothesis,  $\Lie C_G(S)=\cc_\gg(S)$,
so $\Lie C_G(S)$ is the trivial $S$-isotypic summand.  Hence
there is a unique $S$-stable complement to $\Lie C_G(S)$ in $\gg$,
namely, the sum $\mm$ of the non-trivial $S$-isotypic summands. The
uniqueness of $\mm$ implies that it is also $C_G(S)$-stable, so it
follows that $(G,C_G(S))$ is a reductive pair.

(b).
Now assume that $S$ is a subgroup of $G$.
Then $S$ is $G$-cr by Theorem \ref{thm:adjcr}, so $N_G(S)$ is
$G$-cr \cite[Thm.~3.14]{BMR}. Moreover, $N_G(S)/SC_G(S)$ is
finite by \cite[Lem.\ 6.8]{martin1}. Since $S^0\subseteq C_G(S)$ by
assumption, $N_G(S)/C_G(S)$ is finite. It follows that $\Lie
N_G(S)=\Lie C_G(S)$.  By the uniqueness of the subspace $\mm$
of $\gg$ from the proof of part (a) above, $\mm$ is also
$N_G(S)$-stable.  Hence $(G,N_G(S))$ is a reductive pair.
\end{proof}

The following is immediate by
Theorem \ref{thm:redpairs}, Corollary~\ref{cor:redpairs} and Proposition~\ref{prop:C_G(S)redpair}(a).
Note that it applies in particular when $S$ is a torus, so that $C_G(S)$ is an R-Levi subgroup
of $G$ \cite[Cor.~6.10]{BMR}.

\begin{cor}
\label{cor:LGsep}
Let $G$ and $S$ be as in the statement of Proposition \ref{prop:C_G(S)redpair}.
Then every subgroup of $C_G(S)$ which is separable in $G$
is separable in $C_G(S)$ and every subalgebra of $\cc_\gg(S)$
which is separable in $\gg$ is separable in $\cc_\gg(S)$.
%Let $H$ be a subgroup of $C_G(S)$.
%If $H$ is separable in $G$, then $H$ is separable in $C_G(S)$.
\end{cor}

For a regular reductive subgroup $H$ of $G$ the next lemma gives
a useful criterion for $(G,H)$ to be a reductive pair
(cf.~\cite[Ex.~3.33]{BMR}).

\begin{lem}
\label{lem:redpair}
Let $T$ be a maximal torus of $G$ and let $H$ be a reductive subgroup
of $G$ containing $T$. Assume that $\Psi(H)=\Psi(H,T)$ is a closed
subsystem of $\Psi=\Psi(G,T)$. Then $(G,H)$ is a reductive pair.
\end{lem}

\begin{proof}
Let $\mm$ be the sum of the root spaces $\uu_\alpha$ with
$\alpha\notin\Psi(H)$. Then $\gg=\hh\oplus\mm$.
By the conjugacy of the maximal tori in $H^0$,
we have $H=H^0N_H(T)$. Since $N_H(T)$ permutes the root spaces in $\hh$
and also those outside $\hh$, $N_H(T)$ stabilizes $\mm$.
So all we have to show is that $\mm$ is stable
under the $U_\beta$, for $\beta\in\Psi(H)$.

Let $\alpha\in\Psi\setminus\Psi(H)$ and $\beta\in\Psi(H)$.
If $\gamma = \alpha+i\beta$ is a root for some integer $i\ge0$,
then we must have $\gamma\notin\Psi(H)$, as otherwise
$\alpha = \gamma-i\beta\in\Psi(H)$,
because $\Psi(H)$ is closed and symmetric;
see \cite[Ch.~6 Prop.~23(iii)]{bou}.
Now let $u\in U_\beta$.
Then, by \cite[Lem.~5.2]{borel1},
$\Ad u(\uu_\alpha)\subseteq\bigoplus\uu_{\alpha+i\beta}$,
where the sum is over all integers $i\ge0$ such that $\alpha+i\beta$ is a root.
By the above, this sum is contained in $\mm$.
\end{proof}
%Bourbaki, Groupes et Algebres de Lie,
%Ch. 6 par.1 no. 7 Prop. 23(iii).

\begin{rems}
\label{rem:C_G(S)redpair}

(i).
Note that if $S$ is a linearly reductive group acting on $G$ by
automorphisms, then the conditions of Proposition~\ref{prop:C_G(S)redpair}(a)
are satisfied, \cite[Lem.\ 4.1]{rich0}.

(ii).
We can apply Proposition~\ref{prop:C_G(S)redpair}(b) and Theorem~\ref{thm:redpairs} to prove that if $S$ is a linearly
reductive subgroup of $G$ with $S^0$ central in $S$ and $H$ is a
$G$-cr subgroup of $N_G(S)$ such that $H$ is separable in $G$, then
$H$ is $N_G(S)$-cr.  Here the separability condition cannot be
removed; see Example~\ref{exmp:KnotKS}.

(iii).
Note that the hypothesis on $S^0$ in Proposition
\ref{prop:C_G(S)redpair}(b) cannot be removed.  For example, suppose
$p=3$ and let $G= (k^*\times k^*\times k^*)\rtimes (C_3\times C_2)$, where the
cyclic group $C_3$ acts on $k^*\times k^*\times k^*$ by a cyclic permutation of the factors and the cyclic group $C_2= \{1,a\}$ acts on $k^*\times k^*\times k^*$ by $a\cdot (x,y,z)= (x^{-1},y^{-1},z^{-1})$.  Let $S$ be the linearly reductive subgroup
$\Delta C_2$, where $\Delta$ is the diagonal inside $k^*\times k^*\times k^*$.
Then $N_G(S)^0= \Delta$ and $C_3\subseteq N_G(S)$, so $(G,N_G(S))$ is not a
reductive pair because $\Lie \Delta$ does not admit a $C_3$-stable
complement in $\Lie (k^*\times k^*\times k^*)$.

(iv).
We observe that Slodowy's example \cite[I.3(7)]{slodowy} is a
special case of Proposition \ref{prop:C_G(S)redpair}(a), namely when $G
= \GL(V)$ and $C_G(S)$ is a Levi subgroup of $G$.

(v).
In the special case of Proposition \ref{prop:C_G(S)redpair} when
$S$ is linearly reductive and $G$ is connected, Richardson showed
in \cite[Prop.\ 10.1.5]{rich0} that $C_G(S)$ is reductive and, if $S$ is a subgroup of $G$,
$N_G(S)$ is reductive.

(vi).
The converse of Corollary~\ref{cor:LGsep} is false: see
Proposition \ref{prop:LsepGnonsep}.
\end{rems}

We give another application of Proposition \ref{prop:C_G(S)redpair}.

\begin{exmp}
\label{ex:d4g2}
Suppose that $p = 2$.  Let $G$ be simple of type $D_4$ and let $S$ be the
group of order $3$ generated by the triality graph automorphism of
$G$.  Then $K = C_G(S)^0$ is of type $G_2$.  Since $S$ is linearly
reductive, Proposition \ref{prop:C_G(S)redpair}(a) implies that $(G, K)$
is a reductive pair.  In Section \ref{sec:ex}, we construct a subgroup
$H$ of $K$ isomorphic to $S_3$ which is not separable in $K$, see
Proposition \ref{prop:LsepGnonsep}. It follows from Theorem
\ref{thm:redpairs} that $H$ is also non-separable in $G$.  In
addition, by Lemma \ref{lem:HGcr}(a), $H$ is $K$-cr and thus, thanks
to \cite[Cor.\ 3.21]{BMR}, $H$ is also $G$-cr.

This example also gives rise to a non-separable subgroup in the
exceptional group of type $F_4$ as follows.
Let $G'$ denote this group; then,
since $D_4$ is a closed subsystem of the root system of type $F_4$
(the $D_4$-subsystem consists of the long roots in the $F_4$ system),
Lemma \ref{lem:redpair} implies that $(G',G)$ is a reductive pair.
It follows from Theorem \ref{thm:redpairs}
that $H$ is also non-separable in the
group $G'$.
%
%This example also gives rise to a non-separable subgroup in the
%exceptional groups of type $E_6$, $E_7$, and $E_8$
%in characteristic $2$ by means of embedding $D_4$
%as a Levi subgroup and applying \cite[I.3(4)]{slodowy}.
\end{exmp}

\section{The Adjoint Module and Complete Reducibility}
\label{sec:sep1}
%All results in this section are trivial if $p=0$, therefore we will assume throughout that $p>0$.

%Recall Theorem \ref{thm:adjcr}.

In the proof of Theorem \ref{thm:adjcr}, the hypothesis of separability
is used only for a rather coarse dimension-counting argument, so it is natural
to ask whether it can be removed.
This is a more subtle problem than it at first appears.

\begin{comment}
\begin{exmp}[Serre]
\label{ex:serre}
Let $H$ be an exceptional simple group of adjoint type and
suppose $k$ is a field of positive characteristic
$p < 2h-2$, where $h$ denotes the Coxeter number of $H$.
Then, with the possible exception of type $G_2$ in characteristic $5$,
there exists a simple subgroup $K$ of $H$ of type $A_1$ such that
$K$ is $H$-irreducible, but $K$ does not act semisimply on $\hh$ --- that is, $K$ is not $G$-completely reducible for $G:= \GL(\hh)$.
We can choose $p$ and $H$ such that $p$ is good for $H$; then $(G,H)$ is a reductive
pair, since the Killing form on $\hh$ is non-degenerate \cite[\S 5]{rich2}.  It then follows from Theorem \ref{thm:redpairs} that $K$ is separable in $H$, since $K$ is separable in $G$.
%in fact we can arrange
%for $K$ to also be separable in $H$
%(by \cite[Thm.\ 2]{lawthertesterman},
%taking $p>7$ will do).
\end{exmp}
\end{comment}

Our first result shows that
we can remove the separability assumption
from $H$ in Theorem \ref{thm:adjcr}
under extra hypotheses on $H$.

\begin{thm}
\label{thm:allregss}
Suppose that $H$ is a subgroup of $G$ such
that $H$ acts semisimply on ${\mf m}_{\rm ad}$
for every reductive subgroup $M$ of $G$ that contains $H$
and a maximal torus of $G$ (this includes the case $M = G$).
Then $H$ is $G$-completely reducible.
\end{thm}

\begin{proof}
Suppose $P$ is an $R$-parabolic subgroup of $G$ containing $H$,
let $T$ be a maximal torus of $P$
and let $\lambda\in Y(T)$ such that $P=P_\lambda$.  By Remark \ref{rem:P0G0}, we can assume that $P^0$ is proper in $G^0$.  This implies that $\lambda$ is non-central in $G^0$ \cite[Lem.~2.4]{BMR}.

First assume that $H$ acts semisimply on $\gg$ and
that $\zz(\gg)=\{0\}$. We then show
that there exists a reductive subgroup $M$ of $G$
such that $\dim M<\dim G$ and $M$ contains $H$,
$\lambda(k^*)$ and a maximal torus of $G$.
After that we prove the statement of the theorem
using induction on the dimension of $G$.

Put $\nn=\Lie R_u(P)$, $\frakl_\lambda=\Lie L_\lambda $
and $\mf{s}=\Lie\lambda(k^*)$.
Since $L_\lambda$ centralizes $\lambda(k^*)$ and normalizes $R_u(P)$,
$\lambda(k^*)R_u(P)$ is normalized by $L_\lambda$.
Therefore, $\lambda(k^*)R_u(P)$ is a normal subgroup of $P$.
So $\Lie (\lambda(k^*)R_u(P)) = \nn\oplus\mf{s}$
is an $H$-submodule of $\gg$, since $H$ is contained in $P$.
Now $\nn$ is an $H$-submodule of $\nn\oplus\mf{s}$, so,
since $\gg$ is a semisimple $H$-module,
there exists a $1$-dimensional $H$-submodule $\mf{s}_1$ of $\nn\oplus\mf{s}$
which is a direct $H$-complement to $\nn$.

As $L_\lambda$ acts trivially on $\zz(\frakl_\lambda)$,
%by \cite[Lem.~22.2]{borel},
$(\nn\oplus \zz(\mf{l}_\lambda))/\nn$
is a trivial $L_\lambda$-module. By \cite[Prop.~I.3.17]{borel},
we have $\Ad(u)(x)-x\in[\frakl,\nn]\subseteq\mf{n}$ for all
$u\in R_u(P)$ and all $x\in\pp$, so $(\nn\oplus \zz(\mf{l}_\lambda))/\nn$
is also a trivial $R_u(P)$-module and therefore a trivial $P$-module.
So $\mf{s}_1$ must be a trivial $H$-module.
Let $x\in\mf{s}_1$ be non-zero.
Clearly, $H$ fixes the nilpotent and semisimple parts of $x$.
By \cite[Thm.~III.10.6]{borel}, $\nn$ is the set of nilpotent
elements of $\nn\oplus\mf{s}$.
So $x$ has non-zero semisimple part and we may assume that $x$ is semisimple.
Thus we have found a non-zero semisimple element of $\Lie P$ that is fixed by $H$ and
$\lambda(k^*)$.
After conjugating $\lambda$ and $T$ by the same element of
$P$, we may assume that $x\in\Lie T$.

Put $M:=C_{G}(x)$.
This is a reductive subgroup of $G$
which contains $T$, $H$ and $\lambda(k^*)$.
Furthermore, $\dim M < \dim G$,
since otherwise $G^0\subseteq M$ which is impossible,
because $\cc_\gg(G^0) = \zz(\gg)=\{0\}$.
So $M$ has the required properties.

To prove the assertion of the theorem, we pass to the adjoint group.
Let $\pi_G:G\to G_{\rm ad}$ be the natural morphism and let
$\tilde\lambda=\pi_G\circ\lambda\in Y(G_{\rm ad})$.
By our hypothesis, $\pi_G(H)$ acts semisimply on $\gg_{\rm ad}$.
Note that $\tilde\lambda$ is non-trivial.
If we apply the above argument to $G_{\rm ad}$, $\pi_G(H)$ and
$\tilde\lambda$, we
get a reductive subgroup $M'$ of $G_{\rm ad}$ such that
$\dim M'<\dim G_{\rm ad}$ and $M'$ contains $\pi_G(T)$,
$\pi_G(H)$ and ${\tilde\lambda}(k^*)$.
But then $M:=\pi_G^{-1}(M')$ is a reductive subgroup of
$G$ such that $\dim M<\dim G$ and $M$ contains $T$, $H$ and
${\lambda}(k^*)$. Clearly, $M$ and $H$ satisfy the same assumptions as $G$
and $H$.
By the induction hypothesis, $H$ is $M$-cr, so there exists $u\in P_\lambda(M)$
such that $H\subseteq uL_\lambda(M)u^{-1}$.
But $R_u(P_\lambda(M))=R_u(P_\lambda)\cap M$ and similarly for $L_\lambda(M)$.
So $H$ is contained in the Levi subgroup $uL_\lambda u^{-1}$ of $P$.
\end{proof}

\begin{comment}
The following example indicates the sort of problem that can arise
without the assumptions on $H$ made in Theorem \ref{thm:allregss}.
\end{comment}

The following example indicates the sort of problem that can arise
without the assumptions on $H$ made in Theorem \ref{thm:allregss}.  First we need some terminology. If $G$ is
connected and simple of type $A_{n-1}$, then either $\Lie G
\cong\mf{sl}_n$, $\Lie G \cong\mf{pgl}_n$, or $p^2|n$ and $\Lie G$
is of {\it intermediate type}; see \cite[Table~1]{hogeweij}.  In the latter case, $\gg$ is the direct sum of its centre, which is $1$-dimensional, and its derived algebra, which is isomorphic to $\mf{psl}_n:=\mf{sl}_n/(k\cdot{\rm id})$.

\begin{exmp}
\label{ex:gvsgtilde}
Let $p$ be a prime, put $n=p^2$ and let $G$ be the simple algebraic
group of type $A_{n-1}$ whose character group is the lattice that is
strictly between the root lattice and the weight lattice
(the quotient group of the latter two lattices is a cyclic group of order $p^2$).  Then $\gg$ is of intermediate type.
Thus $\gg$ is the direct sum of two simple $G$-modules and is therefore
a semisimple $G$-module.  But the only proper non-zero ideal of $\gg_{\rm ad}\cong\mf{pgl}_n$ is its derived subalgebra, which is of dimension $n^2-2$
(\cite[Table~1]{hogeweij}).  Since every $G$-submodule of $\gg_{\rm ad}$ is an ideal, it follows that $G$ does not act semisimply on $\gg_{\rm ad}$.
\end{exmp}

\begin{comment}
 \begin{exmp}
\label{ex:gvsgtilde}
Let $p$ be a prime, put $n=p^2$ and let $G$ be the simple algebraic
group of type $A_{n-1}$ whose character group is the lattice that is
strictly between the root lattice and the weight lattice
(the quotient group of the latter two lattices is a cyclic group of order $p^2$).
Then $\gg$ is the direct sum of its centre, which is $1$-dimensional,
and its derived
algebra, which is isomorphic to $\mf{psl}_n$ (see \cite[Table~1]{hogeweij}).
Thus $\gg$ is the direct sum of two simple $G$-modules and is therefore
a semisimple $G$-module.
But the only proper non-zero ideal of $\gg_{\rm ad}\cong\mf{pgl}_n$
is its derived subalgebra, which is of dimension $n^2-2$
(\cite[Table~1]{hogeweij}).
Since every $G$-submodule of $\gg_{\rm ad}$ is an ideal,
it follows that $G$ does not act semisimply on $\gg_{\rm ad}$.
\end{exmp}

\end{comment}

\begin{rems}
\begin{comment}
 \begin{itemize}
  \item[(i)] Example \ref{ex:gvsgtilde} shows that it is possible for a subgroup
$H$ to act semisimply on $\gg$,
but not on $\gg_{\rm ad}$.
Thus the argument in the proof of Theorem \ref{thm:allregss}
breaks down in general.
  \item[(ii)] Suppose $G$ is semisimple (so $Z(G)$ is finite),
and $H$ is a subgroup of $G$ which acts semisimply on $\gg$, but
$H$ is not separable in $G$.
Observe that
$H$ might be ``trivially'' non-separable in $G$, in the sense that $\gg$
has non-zero centre so that $G$ is not even separable in itself.
One might hope to remove such possibilities by passing to the
adjoint form $G_{\rm ad}$ of $G$, but
Example \ref{ex:gvsgtilde} shows that if we do this, the image
of $H$ in $G_{\rm ad}$ may fail to act
semisimply on $\gg_{\rm ad}$.
 \end{itemize}
\end{comment}
  (i). Example \ref{ex:gvsgtilde} shows that it is possible for a subgroup
$H$ to act semisimply on $\gg$,
but not on $\gg_{\rm ad}$.
Thus the argument in the proof of Theorem \ref{thm:allregss}
cannot be used to extend Theorem \ref{thm:adjcr} to the non-separable case.%\smallskip\\

  (ii). Suppose $G$ is semisimple (so $Z(G)$ is finite),
and $H$ is a subgroup of $G$ which acts semisimply on $\gg$, but
$H$ is not separable in $G$.
Observe that
$H$ might be ``trivially'' non-separable in $G$, in the sense that $\gg$
has non-zero centre so that $G$ is not even separable in itself.
One might hope to deal with such possibilities by passing to the
adjoint form $G_{\rm ad}$ of $G$, but
Example \ref{ex:gvsgtilde} shows that if we do this, the image
of $H$ in $G_{\rm ad}$ may fail to act
semisimply on $\gg_{\rm ad}$.
\end{rems}

Our next result shows that we can also remove the separability
assumption in Theorem \ref{thm:adjcr} by strengthening the conditions
on $G$, rather than on $H$ (as in Theorem \ref{thm:allregss}).  In
contrast to Theorem \ref{thm:verygood}, this next result does not
impose any characteristic restrictions stemming from simple factors of
$G$ of type $A_n$.

\begin{thm}
\label{thm:adjsc}
Assume that $G$ is connected, $p$ is good for $G$, and no simple
factor of type $A_n$ of the derived group $DG$ of $G$ has Lie algebra
of intermediate type.  Let $H$ be a subgroup of $G$ which acts
semisimply on $\Lie DG$.  Then $H$ is $G$-completely reducible.
\end{thm}

\begin{proof}
Since $Z(G)^0$ acts trivially on $\gg$ and since it is contained in any Levi
subgroup of any parabolic subgroup of $G$, we may assume that $Z(G)^0\subseteq H$.
Then $H=Z(G)^0(H\cap DG)$.
%Applying \cite[Lem.\ 2.12(ii)(b)]{BMR} to the isogeny
%$Z(G)^0\times DG\to G$ and then (i) of the same lemma to $Z(G)^0\times DG$,
Applying \cite[Lem.\ 2.12]{BMR} to the isogeny $Z(G)^0\times DG\to G$ and the projection $Z(G)^0\times DG\ra DG$,
we see that we may replace $G$ by $DG$ and $H$ by $H\cap DG$.

Let $G_1,G_2,\ldots, G_r$ be the simple factors of $G$ and let $\mu : \Pi_i G_i\to G$ be the isogeny given by multiplication.
Denote the projection $G\to G_i$ by $\pi_i$ and put $\gg_i=\Lie G_i$.
Note that $\gg_i$ is $H$-semisimple, since it is a $G$-submodule of $\gg$.
It is easily checked that for each $i$ the set of automorphisms of $\gg_i$ given by $H$
is the same as that given by $\pi_i(\mu^{-1}(H))$. So $\pi_i(\mu^{-1}(H))$ acts
semisimply on $\gg_i$ for each $i$. Thanks to \cite[Lem.\ 2.12]{BMR},
we may now assume that $G$ is simple.

First assume that $G$ is not of type $A$.
Then $H$ is $G$-cr by Theorem \ref{thm:verygood}.

Now assume that $G$ is of type $A_{n-1}$.
%If $p\not| n$ then the desired result follows from Theorem \ref{thm:verygood}, so assume that $p|n$.
Since $\gg$ is not of
intermediate type, the isogenies $\varphi:\SL_n\to G$ and
$\psi:G\to\PGL_n$ cannot both be inseparable. So $\varphi^{-1}(H)$
acts semisimply on $\mf{sl}_n$ or $\psi(H)$ acts semisimply on
$\mf{pgl}_n$. By \cite[Lem.\ 2.12(ii)(b)]{BMR}, we may assume that
$G\cong\SL_n$ or that $G\cong\PGL_n$.

First assume that $G=\SL_n$.  We have $G_{\rm ad}=\PGL_n$ and
$\gg_{\rm ad}=\mf{pgl}_n$.  The trace form on $\gl_n$ is
non-degenerate and induces a non-degenerate $\GL_n$-invariant pairing between
$\mf{sl}_n$ and $\mf{pgl}_n$.  This shows that
$\mf{pgl}_n\cong\mf{sl}_n^*$ as $\GL_n$-modules and therefore that $H$
acts semisimply on $\gg_{\rm ad}$.  More generally, thanks to
%the algorithm of Borel and de Siebenthal
\cite{BoSe}, if $M$ is any reductive subgroup of $G$ containing a
maximal torus of $G$, then $M^0$ is a Levi subgroup of $G$ (see also
\cite[Ex.\ Ch.\ VI \S 4.4]{bou}).  Regarding $G$ as a subgroup of
$\GL_n$, we can write $M^0= (\prod_{i=1}^r\GL_{n_i})\cap\SL_n$ and
$\mm= (\oplus_{i=1}^r\gl_{n_i})\cap\mf{sl}_n$ for some $n_i$.  Since
the restriction of the trace form of $\gl_n$ to
$\oplus_{i=1}^r\gl_{n_i}$ is non-degenerate --- it is the direct sum
of the trace forms of the $\gl_{n_i}$ --- the orthogonal complement of
$\mm $ in $\oplus_{i=1}^r\gl_{n_i}$ is 1-dimensional and is therefore
equal to $k\cdot {\rm id}_n$. So $\mm^*\cong
(\oplus_{i=1}^r\gl_{n_i})/(k\cdot{\rm id}_n)$ as $M$-modules.  Now
$\mm_{\rm ad}=\oplus_{i=1}^r\mf{pgl}_{n_i}$, which is isomorphic to a
quotient of $\mm^*=(\oplus_{i=1}^r\gl_{n_i})/(k\cdot {\rm id}_n)$.
\begin{comment}
\hfill\break BEGIN COMMENT\\ If $L =
(\GL_{n_1}\times\cdots\times\GL_{n_r})\cap\SL_n$ is a Levi of $\SL_n$,
then $\Lie L = \oplus_{i=1}^r\gl_{n_i}\cap\mf{sl}_n$. Since the
restriction of the trace form of $\gl_n$ to $\oplus_{i=1}^r\gl_{n_i}$
is non-degenerate (it is the direct sum of the trace forms of the
$\gl_{n_i}$), the orthogonal of $\Lie L $ in $\oplus_{i=1}^r\gl_{n_i}$
is one dimensional and therefore equal to $k{\rm id}_n$. So $\Lie
L^*\cong (\oplus_{i=1}^r\gl_{n_i})/k{\rm id}_n$ as
$N_{\GL_n}(\GL_{n_1}\times\cdots\times\GL_{n_r})$-modules.  Now $\Lie
\Ad(L)=\oplus_{i=1}^r\mf{pgl}_{n_i}$ which is clearly isomorphic to a
quotient of $(\oplus_{i=1}^r\gl_{n_i})/k{\rm id}_n$.\\ END COMMENT\\
\end{comment}
Thus if $H$ acts semisimply on $\gg$, then $H$ acts semisimply on
$\mm_{\rm ad}$ for any reductive subgroup $M$ which contains
$H$ and a maximal torus of $G$.  The result now follows from Theorem
\ref{thm:allregss}.

Finally, assume that $G$ is isomorphic to $\PGL_n$.  Let
$\varphi:\SL_n\to\PGL_n$ be the canonical projection.  Since
$\mf{pgl}_n\cong\mf{sl}_n^*$ as $\GL_n$-modules, $\varphi^{-1}(H)$
acts semisimply on $\mf{sl}_n$.  The desired result follows from the
previous case and \cite[Lem.\ 2.12(ii)(b)]{BMR}.
\end{proof}

We record two important special cases of Theorem \ref{thm:adjsc}.

\begin{cor}\label{cor:adjc}
Assume that $G$ is connected, $p$ is good for $G$, and $DG$ is either
adjoint or simply connected.  Let $H$ be a subgroup of $G$
which acts semisimply on $\Lie DG$. Then $H$ is $G$-completely
reducible.
\end{cor}

\begin{rems}
(i). Note that Theorem~\ref{thm:adjsc} applies in cases when Theorem~\ref{thm:adjcr}
does not. For example, let $G = \SL_2(k)$, where $k$
has characteristic $2$, let $T$ be a maximal torus of $G$, and set $H
= N_G(T)$.  Then $\gg$ is semisimple as an $H$-module, so Theorem
\ref{thm:adjsc} applies, but, as $H$ is not separable in $G$ (cf.\
\cite[Ex.\ 3.4(b)]{martin3}), Theorem \ref{thm:adjcr} does not apply.
We can also take $G=\GL_2(k)$, $k$ of characteristic $2$, and $H =
N_G(T)$, where $T$ is a maximal torus of $G$. Then $\Lie DG =
\mf{sl}_2$ is a semisimple $H$-module, but $\gg = \gl_2$ is not
semisimple as an $H$-module.

(ii). Let $H$ be a group and let $\rho:H\to\GL(V)$ be a finite-dimensional representation of $H$.  We have $V\otimes V^*\cong\gl(V)$
as $\GL(V)$-modules and therefore also as $H$-modules. Furthermore, $V$
is $H$-semisimple if and only if $\rho(H)$ is $\GL(V)$-cr. So, by
Theorem~\ref{thm:adjsc} (or Theorem \ref{thm:adjcr}), $V$ is
$H$-semisimple if $V\otimes V^*$ is $H$-semisimple. This result is a
special case of a theorem of Serre, see \cite[Thm.~3.3]{serre0}.

(iii). We do not know whether the assumptions other than $\gg$
being $H$-semisimple are necessary in Theorem~\ref{thm:adjsc}.
To show that the assumption on the simple type $A$ factors
can be removed, it would
suffice to show that, under the assumption that $p^2|n$, a subgroup $H$ of $\GL_n$ which acts
semisimply on $\mathfrak{psl}_n$
is $\GL_n$-cr --- that is, acts semisimply on the natural
module $k^n$.
%R: Assume we have $G$ simple of type $A_n$ with Lie algebra of intermediate type and
%$H\subseteq G$ acting semisimply on $\gg$. Moving this group to \SL_n\subseteq\GL_n we
%obtain a subgroup of  group of \GL_n with the properties above. This group would then be
%\GL_n-cr which should imply that H is G-cr.
\begin{comment}
Conversely, to show that the assumption on the simple type $A$
factors cannot be removed, it
suffices to show that the above assertion is false when $p^2|n$.
\end{comment}

(iv). It is easy to prove that $k^n$ is
$H$-semisimple if $\mf{pgl}_n$ (equivalently, $\mf{sl}_n$) is $H$-semisimple.
The arguments are straightforward modifications of the arguments at the
end of the proof of Weyl's theorem in \cite[Thm.~6.3]{Hum2} and work for an arbitrary field $k$.
So one can prove Theorem~\ref{thm:adjsc} without using Theorem~\ref{thm:allregss}.
\end{rems}

We finish this section with an example which illustrates that the
converse of Theorem \ref{thm:adjcr} is false.

\begin{exmp}
\label{exmp:Serre}
In \cite[Cor.\ 5.5]{serre2}, Serre
shows that $\gg$ is a semisimple $H$-module if and only if
$H$ is $G$-completely reducible, under the assumption that
$p > 2h-2$, where $h$ is the Coxeter number of $G$.
Theorem \ref{thm:verygood} improves the bound for the forward implication
of this equivalence.
However, Serre's bound is sharp for the reverse implication,
as the following example (due to some unpublished calculations of Serre)
shows, see also \cite[Rem.~3.43(iii)]{BMR}.

Let $H$ be an exceptional
simple group of adjoint type and suppose
$p \leq 2h-2$, where $h$ denotes the Coxeter number of $H$.
Then, with the possible exception of type $G_2$ in characteristic $5$,
Serre's construction yields a simple subgroup
$K$ of $H$ of type $A_1$ such that
$K$ is $H$-irreducible, but $K$ does not act
semisimply on $\hh$ --- that is, $K$ is not
$G$-completely reducible, where $G:=\GL(\hh)$.

This example also illustrates that the converse of the second
assertion of Theorem \ref{thm:redpairs} is not true in general,
even when $K$ is separable in $G$ and $(G,H)$
is a reductive pair.
For we can choose $p$ and $H$ above such that $p$ is
good for $H$; then $(\GL(\hh),H)$ is a reductive
pair (cf.\ \cite[Ex.\ 3.37]{BMR}), since the Killing form on $\hh$ is
non-degenerate \cite[\S 5]{rich2}.  It
then follows from Theorem \ref{thm:redpairs}
that $K$ is separable in $H$, since $K$ is separable in $\GL(\hh)$.
\end{exmp}

%\begin{exmp}
%Suppose there exists a subgroup $H$ of $G$ such that
%$H$ acts semisimply on $\gg$, but $H$ is not $G$-completely reducible
%(i.e.\ a counterexample to the general case of this section).
%%Section \ref{sec:sep1}).
%Then $(G,H)$ is a reductive pair and $K:=H$ is $H$-completely reducible, but
%not $G$-completely reducible, so we obtain another counterexample to the converse
%of Theorem \ref{thm:redpairs}.
%%(It follows, for two different reasons, that $H$ cannot be
%%separable in $G$).
%\end{exmp}

%%%%%%%%%%%%%%%%%%%%%%%%%%%%%%%%%%%%%%%%%%%%%%%%%%%%%%%%%%%%%%%%%%%%%%%%%%%%%

\section{Closed Orbits and Separability}
\label{sec:orbitthm}

%\subsection{A Result of Richardson}
%\label{sub:Richardson}
In this section we generalize two results of R.W.~Richardson
\cite[Thm.\ A, Thm.\ C]{rich0},
which we are then able to apply to $G$-complete
reducibility.
%Note that our proofs are in the spirit of Richardson's
%original arguments.
%We split our version into two parts, since one direction holds
%in greater generality than the other.  Note that if $S$ is linearly reductive %then we can replace $G'$ by the reductive group $\langle S,G\rangle$, so %Theorem \ref{thm:onedirection} yields \cite[Thm.\ C]{rich0} as a special case.
%{\tt R:In Theorem~\ref{thm:richC} b) one could replace "the adjoint
%representation of $S$ on $\gg$ is semisimple" by "all irreducible
%$S$-constituents of $\gg/\cc_\gg(S)$ are nontrivial".}\\
%
%\begin{comment}%
%BEGIN COMMENT\\
We require two preliminary results,
which allow us to relax the hypotheses in the
original theorems.
Armed with these lemmas, it is quite straightforward to adapt
Richardson's proofs to our setting.

Our first result is an analogue of \cite[Prop.\ 6.1]{rich0},
which Richardson proves for $S$ linearly reductive
and $G$ connected, using non-commutative cohomology of
algebraic groups.
It should be noted that, although our result is more general
than Richardson's, in order to apply it to linearly reductive groups
one needs to know that they are always completely reducible in
any ambient reductive group (\cite[Lem.\ 2.6]{BMR}),
and the proof of this depends on the argument that Richardson uses.

\begin{lem}
\label{lem:Sstable}
Let $G$ be a normal subgroup of the reductive group $G'$ and
suppose that $S$ is a $G'$-completely reducible subgroup of $G'$.
Then for every R-parabolic subgroup $P$ of $G$ normalized by $S$,
there exists an R-Levi subgroup of $P$ normalized by $S$.
\end{lem}

\begin{proof}
Note that since $G$ is normal in $G'$, $G$ is reductive,
and $G^0$ is a connected reductive normal subgroup of $G'$.
We first show that it is enough to prove the result when $G$ is connected.
If $P$ is an R-parabolic subgroup of $G$ normalized by $S$,
then $P^0$ is a parabolic subgroup of $G^0$ normalized by $S$.
If $L$ is a Levi subgroup of $P^0$ normalized by $S$, then
$L = M^0$ for some R-Levi subgroup $M$ of $P$ \cite[Cor.~6.8]{BMR}.
Let $T$ be a maximal torus of $M$; then $T \subseteq L$.
Since $S$ normalizes $P$, it follows that $sMs\inverse$ is
an R-Levi subgroup of $P$ for any $s \in S$.
Moreover, since $S$ normalizes $L$, we see that
$T \subseteq L= sLs\inverse\subseteq sMs\inverse$.
This shows that $M = sMs\inverse$, for every $s \in S$,
by \cite[Cor.\ 6.5]{BMR},
and therefore $S$ normalizes $M$.
Thus, if the result is true for $G^0$, then it is true for $G$.

Now suppose $G$ is connected, and let $H = C_{G'}(G)^0$.
Since $G$ is normal in $G'$, we have $G'^0 = GH$, by
\cite[Lem.\ 6.8]{martin1}.  Let $P$ be a parabolic
subgroup of $G$ normalized by $S$.  Then $P\times H$
is a parabolic subgroup of $G\times H$, and it follows
from \cite[Prop.~IV.11.14(1)]{borel} applied to the
multiplication map $f\colon G\times H\ra G'^0$ that
$PH$ is a parabolic subgroup of $G'^0$.  As $S$
normalizes $P$, $S$ normalizes $PH$, so $S\subseteq N_{G'}(PH)$,
which is an R-parabolic subgroup of $G'$ by \cite[Prop.\ 5.4(a)]{martin1}.
Now $S$ is $G'$-cr by hypothesis, so $S$ is contained in some R-Levi
subgroup $M$ of $N_{G'}(PH)$, so $S$ normalizes $M^0$, which is a Levi
subgroup of $N_{G'}(PH)^0=PH$.  By \cite[Lem.~6.15(ii)]{BMR}, $f^{-1}(PH)$
is a parabolic subgroup of $G\times H$ with Levi subgroup $f^{-1}(M^0)$.

If $(g,h)\in f^{-1}(PH)$, then $gh=g'h'$ for some $g'\in P$
and some $h'\in H$,
so $(g')^{-1}g= h'h^{-1}\in G\cap H\subseteq Z(G)\subseteq P$, so $g\in P$.
It follows that $f^{-1}(PH)=P\times H$.  Hence $f^{-1}(M^0)= L\times H$
for some Levi subgroup $L$ of $P$ and hence $M^0=LH$.
A similar calculation shows that $L=M^0\cap G$, so $S$
normalizes $L$ and we are done.
\end{proof}

\begin{lem}
\label{lem:richC}
Let $S$ be a group and let $f : V\to W$ be a surjective homomorphism of
finite-dimensional $S$-modules over some field. If $V/V^S$ has no trivial
composition factor, then $f(V^S)=W^S$.
\end{lem}

\begin{proof}
Since $f$ is surjective, it induces a surjective homomorphism
$V/V^S\to W/f(V^S)$ of $S$-modules. Since $V/V^S$ has no trivial
composition factor, neither does $W/f(V^S)$. In
particular, $W/f(V^S)$ does not have any non-zero $S$-fixed points.
So $W^S$ must be equal to $f(V^S)$.
\end{proof}

\begin{rem}
\label{rem:richC}
If $V$ is $S$-semisimple ---
in particular, if $S$ is a linearly reductive algebraic group --- then  $V/V^S$ has no
trivial composition factor, so Lemma \ref{lem:richC} applies.
%Recall that a constituent of $V$ is a
%module that is isomorphic to a quotient of a submodule of $V$. The
%irreducible constituents of $V$ are the quotients of a Jordan-H\"older
%series for $V$.
\end{rem}

We can now state our generalization of Richardson's two results
\cite[Thm.\ A, Thm.\ C]{rich0}.
Observe that if $S$ is a subgroup of $G$, then
the condition on $S$ in Theorem \ref{thm:richC}(b)(ii)
is just that of separability of $S$ in $G$.

\begin{thm}
\label{thm:richC}
Let $G$ be a normal subgroup of the reductive group $G'$.  Let $S$ be
a subgroup of $G'$ such that $G'= GS$.  Let $G'$ act on a
variety $X$ and let $x\in X^S$.  Let $\mathcal O$ denote the unique
closed $G$-orbit in the closure of $G\cdot x$.
\begin{itemize}
\item[(a)] Suppose that $X$ is affine and $S$ is $G'$-completely reducible.  Then
$C_G(S)$ is $G$-completely reducible and there exists $\lambda\in
Y(C_G(S))$ such that $\underset{a\to 0}{\lim}\, \lambda(a)\cdot x$
exists and is a point of $\mathcal O$.  In particular, $G\cdot x$ is
closed if $C_G(S)\cdot x$ is closed.
\item[(b)] Suppose that
\begin{itemize}
\item[(i)] $\gg/\cc_\gg(S)$ does not have any trivial $S$-composition
factors;
%the adjoint representation of $S$ on $\gg$ is semisimple;
%\item[(i)] $\cc_\gg(S)$ has an $S$-stable complement in $\gg$; and
\item[(ii)] $\cc_\gg(S)= \Lie C_G(S)$.
\end{itemize}
Then each
%$S$ is $G'$-completely reducible and
irreducible component of $G\cdot x \cap X^S$ is a single
$C_G(S)^0$-orbit. In particular, $C_G(S)\cdot x$ is closed if $G\cdot
x$ is closed.
 \end{itemize}
\end{thm}

\begin{proof}
(a). Suppose that $S$ is $G'$-cr. Then $C_G(S)$ is $G$-cr, by the same
arguments as in the proof of Proposition~\ref{prop:C_G(S)redpair}(a).
If $G\cdot x$ is closed, then we can take $\lambda$ to be the zero
cocharacter, so assume that $G\cdot x$ is not closed.  Then $G^0\cdot
x$ is not closed.  Moreover, the unique closed $G$-orbit $\mathcal{O}$
in the closure of $G\cdot x$ contains the unique closed $G^0$-orbit in
the closure of $G^0\cdot x$.  Thus, since $Y(G) = Y(G^0)$ and
$Y(C_G(S)) = Y(C_{G^0}(S))$, we may replace $G$ by $G^0$ and assume
that $G$ is connected.

We can now follow the first part of Richardson's original proof in
\cite[Sec.~8]{rich0} word for word to give the result, simply
replacing references to \cite[Prop.\ 6.1]{rich0} with our Lemma
\ref{lem:Sstable}.

(b). Now assume that $\gg/\cc_\gg(S)$ does not have any trivial
$S$-composition factors and that $\cc_\gg(S)= \Lie C_G(S)$.  We can
write $G\cdot x$ as a finite union ${\mathcal
O}_1\cup\cdots\cup{\mathcal O}_n$ of $G^0$-orbits. Each ${\mathcal
O}_i$ is a $G$-translate of $G^0\cdot x$, from which we deduce that
the ${\mathcal O}_i$ are the irreducible components of $G\cdot x$. We
may assume that ${\mathcal O}_1,\ldots,{\mathcal O}_r$ $(r\le n)$, are
the $G^0$-orbits in $G\cdot x$ that meet $X^S$. Then $G\cdot x\cap
X^S=({\mathcal O}_1\cap X^S)\cup\cdots\cup({\mathcal O}_r\cap X^S)$, a
disjoint union of closed subsets of $G\cdot x$. This shows that, after
replacing $G'$ by $G^0S$, we may assume that $G$ is connected.
%R:The point is that each irr. component of G\cdot x\cap X^S lives in some {\mathcal O}_i\cap X^S.
%I need the connectedness for the argument at the end of the proof of \cite[Thm.\ A]{rich0}.

Note that $G\cdot x$ is $S$-stable. Put $Y=(G\cdot x)^S = G\cdot x\cap X^S$.
By the argument at the end of the proof of \cite[Thm.\ A]{rich0},
we may assume that the orbit map $G\to G\cdot x$ of $x$ is separable:
Let $y\in Y$.
Then the orbit map of $y$ is separable: that is, its differential
$\psi:\gg\to T_y(G\cdot y)$
is surjective. Since the orbit map of $y$ is $S$-equivariant,
$\psi$ is $S$-equivariant.
By Lemma~\ref{lem:richC}, we have $\psi(\cc_\gg(S))=(T_y(G\cdot y))^S$.
Clearly, $T_y(Y)\subseteq (T_y(G\cdot x))^S$,
so $T_y(Y)\subseteq\psi(\cc_\gg(S))$. Since $\cc_\gg(S)=\Lie C_G(S)$,
the Tangent Space Lemma \cite[Lem.~3.1]{rich0} gives the first assertion of (b).

It now follows that $C_G(S)^0\cdot x$, and therefore $C_G(S)\cdot x$,
is closed in $G\cdot x$.
So if $G\cdot x$ is closed, then $C_G(S)\cdot x$ is closed.
\end{proof}

\begin{rems}\label{rems:richC}
  (i). Assume that $S$ is a separable subgroup of
  $G$ which acts semisimply on $\gg$.
  Then $S$ is $G$-cr by Theorem \ref{thm:adjcr}.
  So parts (a) and (b) of Theorem~\ref{thm:richC} both apply.
  Thus in this case $C_G(S)\cdot x$ is closed if and
  only if $G\cdot x$ is closed.

  (ii). Suppose that $S$ is linearly reductive.  Then conditions (i)
and (ii) of Theorem \ref{thm:richC}(b) are satisfied (see
\cite[Lem.~4.1]{rich0}), and $S$ is $G'$-cr \cite[Lem.~2.6]{BMR};
the special case of
Theorem \ref{thm:richC} for $S$ linearly reductive is precisely
Richardson's original theorem \cite[Thm.\ C]{rich0} (clearly, the
hypothesis in Theorem \ref{thm:richC} that $G'=GS$ is harmless).

(iii). It follows from the proof of Theorem \ref{thm:richC}(b)
that \cite[Thm.~A]{rich0} holds with the hypothesis that $S$
acts semisimply on $\gg$ replaced by the weaker hypothesis of
Theorem \ref{thm:richC}(b)(ii).

(iv). The result in Theorem \ref{thm:richC}(a) that the $G$-orbit
is closed if the $C_G(S)$-orbit is closed is also a slight generalization of part of
\cite[Thm.\ 4.4]{bate}.

(v). Theorem \ref{thm:richC}(b) does not hold for an arbitrary
$G'$-completely reducible subgroup $S$.  In fact it does not even hold
when $S$ is a $G$-completely reducible subgroup of $G$: in
\cite[Prop.~3.9]{BMR2} it is shown that if $S$ is a $G$-completely
reducible subgroup of $G$ and $K$ is a subgroup of $C_G(S)$, then $K$
is $C_G(S)$-completely reducible if and only if $KS$ is $G$-completely
reducible.  Now \cite[Examples 5.1, 5.3, 5.5]{BMR2} give instances of
commuting $G$-completely reducible subgroups $K$ and $S$ such that
$KS$ is not $G$-completely reducible, whence $K$ is not
$C_G(S)$-completely reducible.  If the field $k$ is large enough (cf.\
\cite[Rem.~2.9]{BMR}), we can pick an $n$-tuple $(k_1, \ldots, k_n) \in
G^n$ topologically generating $K$ for some $n$.  Then by Theorem
\ref{thm:gcrcrit}, $G\cdot(k_1, \ldots, k_n)$ is closed in $G^n$, but
$C_G(S)\cdot(k_1, \ldots, k_n)$ is not.
\end{rems}

\vfill
\eject

\begin{comment}
{\tt R: Note that, because of the lemma on ascending sequences of
finite subgroups in the relative preprint, Corollary \ref{cor:SnormalizesK}
is valid without the assumption that $S$ is reductive
and Corollary \ref{cor:richC} is valid without the assumption that $S$ is
reductive in case char $=p>0$.}
\end{comment}

\begin{cor}
\label{cor:richC}
Let $G$, $G'$, $S$, $X$ and $x \in X^S$ be as in
Theorem~\ref{thm:richC}, and suppose $X$ is affine. Suppose that
\begin{itemize}
%\item[(i)] $S$ is reductive;
\item[(i)] the adjoint representation of $S$ on $\gg$ is semisimple;
\item[(ii)] $\cc_\gg(S)= \Lie C_G(S)$.
\end{itemize}
Then $S$ is reductive, $C_G(S)$ is $G$-completely reducible and $G \cdot x$ is closed if
and only if $C_G(S) \cdot x$ is closed.
\end{cor}

\begin{proof}
 Since $GS=G'$, the quotient $S/(S\cap G^0)$ is isomorphic to a finite-index subgroup of $G'/G^0$, and hence is reductive.  Now $S\cap G^0$ is normal in $S$, so $S\cap G^0$ acts semisimply on $\gg$ by Clifford's Theorem.  The kernel $N$ of the $(S\cap G^0)$-action on $\gg$ is a subgroup of the diagonalizable group $Z(G^0)$ and is therefore reductive. Since $N$ and $(S\cap G^0)/N$ are reductive, $S\cap G^0$ is reductive, and it follows that $S$ is reductive.

Suppose $\Char k=0$.  Note that, in this case, properties
%(ii) and (iii)
(i) and (ii) and the $G'$-complete reducibility of $S$ follow from the fact that $S$ is reductive.  Therefore, the hypotheses of
Theorem~\ref{thm:richC}(a) and (b) are satisfied.

Now suppose $\Char k=p >0$.
%Note that $G\rtimes S$ is reductive, since $(G\rtimes S)^0=G^0\rtimes S^0$
%is reductive by \cite[Cor.\ 14.11]{borel}. Therefore we may assume that $G'=G\rtimes S$.
%R: the point is that any $G'$-variety is also a $G\rtimes S$-variety by means of the
%homomorphism $G\rtimes S\to G'$
By the same arguments as in the proof of
Proposition~\ref{prop:C_G(S)redpair}(a),
we may assume that $S$ is finite and deduce that then $S$ is $G'$-cr.
Thus, the hypotheses of Theorem~\ref{thm:richC}(a) and (b) are again
satisfied.
\end{proof}

Now we translate our results in terms of $G$-complete reducibility.

% A new version I decided not to use.
\begin{comment}
\begin{prop}
\label{prop:SnormalizesK}
Let $G$, $G'$ and $S$ be as in Theorem \ref{thm:richC}.  Let $K\subseteq H$ be
subgroups of $G$ such that $S$ normalizes $K$ and $H^0 = C_G(S)^0$.
\begin{itemize}
\item[(a)] Suppose that $S$ is $G'$-completely reducible.  Then
$H$ is reductive and $K$ is $G$-completely reducible if it is $H$-completely reducible.
\item[(b)] Suppose that
\begin{itemize}
\item[(i)] $\gg/\cc_\gg(S)$ does not have any trivial
$S$-composition factors;
\item[(ii)] $\cc_\gg(S)= \Lie C_G(S)$;
\item[(iii)] $C_G(S)$ is reductive;
\item[(iv)] $S$ is reductive.
\end{itemize}
Then $K$ is $H$-completely reducible if it is $G$-completely reducible.
\end{itemize}
\end{prop}
\end{comment}

\begin{prop}
\label{prop:SnormalizesK}
Let $G$ be a normal subgroup of the reductive group $G'$. Let $S$ be
a subgroup of $G'$ such that $G'= GS$. Let $K\subseteq H$ be
subgroups of $G$ such that $H^0 = C_G(S)^0$ and $S$ normalizes $K$.
\begin{itemize}
\item[(a)] Suppose that $S$ is $G'$-completely reducible.  Then
$H$ is reductive and $K$ is $G$-completely reducible if it is $H$-completely reducible.
\item[(b)] Suppose that
\begin{itemize}
\item[(i)] $\gg/\cc_\gg(S)$ does not have any trivial
$S$-composition factors;
\item[(ii)] $\cc_\gg(S)= \Lie C_G(S)$;
\item[(iii)] $C_G(S)$ is reductive.
\end{itemize}
Then $H$ is reductive and $K$ is $H$-completely reducible if it is $G$-completely reducible.
\end{itemize}
\end{prop}

\begin{proof}
In case (a), $C_G(S)$ is reductive by Theorem~\ref{thm:richC}(a) and in
case (b), this is true by assumption. Since $H^0 = C_G(S)^0$, it follows
that $H$ is reductive.

Clearly, we can now assume that $\Char k=p>0$.
By the argument in the proof of
Proposition~\ref{prop:C_G(S)redpair}(a),
we may also assume that $S$ is finite.  As
in \cite[Lem.~2.10]{BMR}, we can replace $K$ with a subgroup $K'$ that
is topologically generated by some $k_1, \ldots, k_n$ with the property
that for any $\lambda\in Y(G)$, we have $K\subseteq P_\lambda$ if and
only if $K'\subseteq P_\lambda$, and $K\subseteq L_\lambda$ if and
only if $K'\subseteq L_\lambda$.  Thus $K$ is $G$-cr (respectively
$H$-cr) if and only if $K'$ is $G$-cr (respectively $H$-cr).

Since $S$ is finite, we may assume, by replacing $(k_1, \ldots, k_n)$
with a larger tuple if necessary, that $S$
permutes the $k_i$ and therefore that $S$ also normalizes $K'$. Since
$H^0 = C_G(S)^0$, we have that $H \cdot (k_1, \ldots, k_n)$ is closed
if and only if $C_G(S) \cdot (k_1, \ldots, k_n)$ is closed.

Let $G'$ act on $G^n$ by simultaneous conjugation.
The symmetric group $S_n$ acts naturally on $G^n$, and the $G'$-action
commutes with this action.  Set $X = G^n/S_n$ and let $\pi: G^n \to X$
be the natural map; the fibres of $\pi$ are precisely the $S_n$-orbits (see \cite[Sec.~2]{BaRi}, for example).  For any subgroup $M$ of $G'$ and any $(g_1 ,
\ldots, g_n) \in G^n$, we have that $M\cdot (g_1 , \ldots, g_n)$ is closed
in $G^n$ if and only if $M\cdot \pi((g_1 , \ldots, g_n))$ is closed in
$X$. Put $x=\pi((k_1, \ldots, k_n))$. Then $x\in X^S$.
Assertions (a) and (b) now follow from
Theorem~\ref{thm:gcrcrit} and from Theorem~\ref{thm:richC}(a) and (b),
respectively.
\end{proof}

\begin{comment}
\begin{rem}
 The hypothesis that $S$ is reductive in Proposition \ref{prop:SnormalizesK} is unnecessary.  This is clear in characteristic 0 --- all that is needed is that $C_G(S)$ is reductive.  In characteristic $p>0$, it can be shown that {\bf any} algebraic group $S$ admits an ascending chain of finite subgroups $S_1\subseteq S_2\subseteq\cdots$
such that $\bigcup_{i\ge1}S_i$ is Zariski dense in $S$, so the argument of Proposition \ref{prop:C_G(S)redpair}(a) holds for non-reductive $S$ as well.
\end{rem}
\end{comment}

The next corollary is a generalization of \cite[Cor.\ 3.21]{BMR}.
Note that the hypotheses on $S$
in
%Corollary~\ref{cor:richC}
Corollary \ref{cor:SnormalizesK} are satisfied if $S$
is linearly reductive.

\begin{cor}
\label{cor:SnormalizesK}
Let $S$ be
%a reductive
an algebraic group acting on $G$ by automorphisms.
Suppose that $S$ acts semisimply on $\gg$ and $\cc_\gg(S)=\Lie C_G(S)$.
Let $K\subseteq H$ be subgroups of $G$ such that $H^0 = C_G(S)^0$ and $S$ normalizes $K$.
Then $H$ is reductive and $K$ is $G$-completely reducible if and
only if $K$ is $H$-completely reducible.
\end{cor}

\begin{proof}
By Proposition~\ref{prop:C_G(S)redpair}(a), $C_G(S)$ is reductive.
So $H$ is also reductive.
We can now assume that $\Char k=p>0$.
By the argument in the proof of
Proposition~\ref{prop:C_G(S)redpair}(a),
%(a),
we can assume that $S$ is finite. Then we put $G'=G\rtimes S$
and obtain, as in that proof, that $S$ is $G'$-cr and that $C_G(S)$ is $G$-cr.
Now the assumptions of Proposition~\ref{prop:SnormalizesK}(a) and (b) are satisfied.
\end{proof}

\begin{comment}
\begin{cor}
\label{cor:SnormalizesK}
Let $S$ be a reductive algebraic group acting on $G$ by automorphisms.
Suppose that $S$ acts semisimply on $\gg$ and $\cc_\gg(S)=\Lie C_G(S)$.
Let $K\subseteq H$ be subgroups of $G$ such that $S$ normalizes $K$
and $H^0 = C_G(S)^0$.
Then $H$ is reductive and $K$ is $G$-completely reducible if and
only if $K$ is $H$-completely reducible.
\end{cor}

\begin{proof}
By Proposition~\ref{prop:C_G(S)redpair}(a), $C_G(S)$ is reductive.
So $H$ is also reductive.
We can now assume that $\Char k=p>0$.
By the argument in the proof of
Proposition~\ref{prop:C_G(S)redpair}(a),
we may assume that $S$ is finite. Then we put $G'=G\rtimes S$
and obtain, as in that proof, that $S$ is $G'$-cr and that $C_G(S)$ is $G$-cr.
Now the assumptions of Proposition~\ref{prop:SnormalizesK}(a) and (b) are satisfied.
\end{proof}
\end{comment}

%%%%%%%%%%%%%%%%%%%%%%%%%%%%%%%%%%%%%%%%%%%%%%%%%%%%%%%%%%%%%%%%%%%%%%%%%%%%%%

\section{Centralizers and Normalizers}
\label{sec:cent+norm}

%{\tt R: This is now the weakest section. We need some examples.}

\begin{comment}
In the previous section we looked at the transfer of
complete reducibility from $G$ to subgroups of $G$ of the form $C_G(S)$,
where $S$ is a group acting on $G$.
We now consider the situation that $S$ is a
subgroup of $G$ acting by inner automorphisms.
This additional requirement on $S$ allows us to
replace $C_G(S)$ by a subgroup $H$ that is closely related to
$C_G(S)$; our next result provides some information about what happens
when we just assume that $H$ contains $C_G(S)^0$, rather than $H^0 =
C_G(S)^0$, as in Proposition \ref{prop:SnormalizesK}.
The result also generalizes \cite[Prop.\ 3.9]{BMR2}, which
deals with the special case when $H=C_G(S)$; in fact, the argument we
use is very similar.
\end{comment}

In this section we continue the theme of Section \ref{sec:orbitthm}, looking at the special case when $S$ is a subgroup of $G$ acting on $G$ by inner automorphisms.  The extra restriction on $S$ allows us to consider subgroups $H$ sitting between $C_G(S)^0$ and $N_G(S)$.  The following result gives a criterion for $K$ to be $H$-cr; it generalizes \cite[Prop.\ 3.9]{BMR2} and the argument we use is very similar.

%MB:- Is this better than a lemma?%
\begin{comment}
\begin{lem}
\label{lem:KN<=>K}
Suppose $S$ is a $G$-completely reducible subgroup of $G$, and suppose
$H$ is a subgroup of $G$ such that $C_G(S)^0 \subseteq H^0 \subseteq
N_G(S)$.  Then $HS$ is reductive and for any subgroup $K \subseteq H
\cap N_G(S)$, $K$ is $H$-completely reducible if and only if $KS$ is
$G$-completely reducible if and only if $KS$ is $N_G(S)$-completely
reducible if and only if $KS$ is $HS$-completely reducible.  If $H$ is
reductive, then $K$ is $H$-completely reducible if and only if $KS$ is
$G$-completely reducible.
\end{lem}
\end{comment}

\begin{comment}
{\tt R: I have now required that $H\subseteq N_G(S)$.
Note: If we only have $H^0\subseteq N_G(S)$, then why would we have $(HS)^0\subseteq N_G(S)$?
What does $HS$ mean in this situation? The closed subgroup generated by $H$ and $S$?}
\end{comment}

\begin{prop}
\label{prop:KN<=>K}
Suppose $S$ is a $G$-completely reducible subgroup of $G$, and suppose
$H$ is a subgroup of $G$ such that $C_G(S)^0 \subseteq H\subseteq
N_G(S)$.  Let $K$ be a subgroup of $H$.
Then:
 \begin{itemize}
  \item[(a)] $HS$ is reductive.
  \item[(b)] The following are equivalent:
 \begin{itemize}
\item[(i)] $KS$ is $G$-completely reducible; % if and only if
\item[(ii)] $KS$ is $HS$-completely reducible; % if and only if
\item[(iii)] $KS$ is $N_G(S)$-completely reducible; % if and only if
\item[(iv)] $KS/S$ is $N_G(S)/S$-completely reducible.
 \end{itemize}
  \item[(c)] Suppose that $H$ is reductive.  Let $\psi$ be the
  canonical projection from $H$ to $H/C_G(S)^0$.  Then $K$ is
  $H$-completely reducible if and only if $KS/S$ is
  $N_G(S)/S$-completely reducible and $\psi(K)$ is
  $H/C_G(S)^0$-completely reducible.
 \end{itemize}
\end{prop}

\begin{proof}
(a). We have
% \begin{equation}
% \label{eqn:norm}
  $C_G(S)^0S^0= N_G(S)^0$,
% \end{equation}
 by \cite[Lem.\ 6.8]{martin1}.
 As $C_G(S)^0\subseteq H$, we have $(HS)^0= N_G(S)^0$.
 Now $N_G(S)$ is reductive by \cite[Prop.\ 3.12]{BMR}, so $HS$ is reductive.

(b). By part (a), $(HS)^0 = N_G(S)^0$, so $HS$ is a finite-index subgroup of $N_G(S)$.
The equivalence of (ii) and (iii) now follows from \cite[Prop.\ 2.12]{BMR2}.
The equivalence of (i) and (iii) follows from \cite[Cor.\ 3.3]{BMR2},
and the equivalence of (iii) and (iv) from \cite[Thm.\ 3.4]{BMR2}.
%Set $N_H(S)= H\cap N_G(S)$.  Then $(HS)^0= N_G(S)^0=
%C_G(S)^0S^0\subseteq N_H(S)S\subseteq HS$, so $N_H(S)S$ is a
%finite-index subgroup of $HS$; in particular, $N_H(S)S$ is reductive.
%We have $N_G(S)^0= C_G(S)^0S^0\subseteq N_H(S)S\subseteq N_G(S)$, so
%$HS$ is a finite-index subgroup of $N_G(S)$.  By \cite[Prop.\
%2.15]{BMR2}, $K$ is $HS$-cr if and only if $K$ is $N_H(S)S$-cr if and
%only if $K$ is $N_G(S)$-cr.  By \cite[Thm.\ 3.4]{BMR2}, $KS/S$ is
%$N_G(S)/S$-cr if and only if $KS$ is $N_G(S)$-cr if and only if $KS$
%is $G$-cr, so (b) follows.

(c). Now suppose that $H$ is reductive.
%We have $C_G(S)^0\subseteq
%N_H(S)^0\subseteq N_G(S)^0= C_G(S)^0S^0$, so $N_H(S)^0= C_G(S)^0(H\cap
%S)^0$.  By a similar argument, $H^0= C_G(S)^0(H\cap S)^0$.  Hence
%$N_H(S)$ is a finite-index subgroup of $H$ and $\psi(N_H(S))$ is a
%finite-index subgroup of $H/C_G(S)^0$.  By \cite[Prop.\ 2.12]{BMR2},
%we can replace $H$ with $N_H(S)$ and $H/C_G(S)^0$ with
%$N_H(S)/C_G(S)^0$, so without loss of generality we assume that
%$H\subseteq N_G(S)$.
The subgroups $H\cap S$ and $(H\cap S)^0$ are normal in $H$.
Let $\pi$ be the canonical projection from $H$ to $H/(H\cap
S)^0$.
We have a homomorphism $\pi\times \psi$ from $H$ to $H/(H\cap
S)^0\times H/C_G(S)^0$.
If $h\in \ker (\pi\times \psi)$ then $h\in
(H\cap S)^0\cap C_G(S)^0\subseteq S^0\cap C_G(S)^0\subseteq Z(S^0)$,
so $\ker (\pi\times \psi)^0$ is a torus, so $\pi\times \psi$ is
non-degenerate.
By \cite[Lem.~2.12(i) and (ii)]{BMR}, $K$ is $H$-cr
if and only if $\pi(K)$ is $H/(H\cap S)^0$-cr and $\psi(K)$ is
$H/C_G(S)^0$-cr.
To finish the proof, it is enough to show that
$\pi(K)$ is $H/(H\cap S)^0$-cr if and only if $KS/S$ is $N_G(S)/S$-cr.
The natural map from $H/(H\cap S)^0\ra HS/S$ is an
isogeny, so $\pi(K)$ is $H/(H\cap S)^0$-cr if and only if $KS/S$ is
$HS/S$-cr, by \cite[Lem.~2.12(ii)]{BMR}.
As $HS/S$ is a finite-index
subgroup of $N_G(S)/S$, $KS/S$ is $HS/S$-cr if and only if $KS/S$ is
$N_G(S)/S$-cr \cite[Prop.\ 2.12]{BMR2}.
This completes the proof.
\end{proof}

\begin{cor}
\label{cor:KS<=>K}
 Let $G$, $H$, $S$ and $K$ be as in Proposition \ref{prop:KN<=>K}.
 Suppose that $(H\cap S)^0$ is a torus.  Then $H$ is reductive and the
 following are equivalent:
 \begin{itemize}
  \item[(i)] $K$ is $H$-completely reducible;
  \item[(ii)] $KS$ is $G$-completely reducible;
  \item[(iii)] $KS$ is $HS$-completely reducible;
  \item[(iv)] $KS$ is $N_G(S)$-completely reducible;
  \item[(v)] $KS/S$ is $N_G(S)/S$-completely reducible.
 \end{itemize}
\end{cor}

\begin{proof}
 Since $N_G(S)^0 = C_G(S)^0S^0$, and $C_G(S)^0 \subseteq H \subseteq N_G(S)$,
 we have $H^0= C_G(S)^0(H\cap S)^0$, so $H$ is
 reductive.  Moreover, $H^0/C_G(S)^0$ is a quotient of $(H\cap S)^0$,
 which is a torus by hypothesis.  Hence $H/C_G(S)^0$ is a finite
 extension of a torus, which implies that any subgroup of $H/C_G(S)^0$
 is $H/C_G(S)^0$-cr.  Thus $K$ is $H$-cr if and only if $KS/S$ is $N_G(S)/S$-cr, by Proposition \ref{prop:KN<=>K}(c).
%Thus one condition in Proposition \ref{prop:KN<=>K}(b)
%is automatic.
% BM: The reference above should be to part (c), not part (b).
 The result now follows from Proposition \ref{prop:KN<=>K}(b).
\end{proof}

\begin{rems}
\label{rem:KvsKN}
(i).
%By \cite[Prop.\ 3.15]{BMR2}, if $N$ is a $G$-cr subgroup of $G$ and $K$ is an $N_G(N)$-cr subgroup of $N_G(N)$,
%then $KN$ is $G$-cr, hence $KN$ is $N_G(N)$-cr by Proposition \ref{prop:KN<=>K}.
%The converse is not true in general: let $N$ be
%$G$-cr but not linearly reductive,
%and let $K$ be a non-$N$-cr subgroup of $N$.
%Then $K$ is not $N_G(N)$-cr by \cite[Prop.~2.8]{BMR2}, but $KN = N$ is $G$-cr.
\begin{comment}
If $H$ is assumed to be reductive, then it follows from Proposition
\ref{prop:KN<=>K} that one set of implications in Corollary
\ref{cor:KS<=>K} above holds without the restriction on $(H\cap S)^0$:
\end{comment}
Suppose $H$ is reductive.  If $K$ is $H$-cr, then Proposition
\ref{prop:KN<=>K}(c) says that $KS/S$ is $N_G(S)/S$-cr, whence $KS$ is
$G$-cr, by Proposition \ref{prop:KN<=>K}(b).  The converse,
however, is not true in this generality; just take $S = H = G$ and $K$
a non-$G$-cr subgroup of $G$.

(ii). Corollary \ref{cor:KS<=>K} holds in particular if $S$ is
linearly reductive, since then the condition that $(H\cap S)^0$ is a
torus is automatic.  However, Example \ref{exmp:KnotKS} shows that
even when $S$ is linearly reductive, the situation for subgroups of $N_G(S)$ is not as
straightforward as for subgroups of $C_G(S)$ (cf.\ Corollary \ref{cor:SnormalizesK}).
%: here we have a linearly reductive subgroup $S$ of a
%(connected) reductive group $G$ and a subgroup $K$ of $N_G
%(S)$ such
%that $K$ is $G$-cr, but $KS$ is not $G$-cr, hence $K$ and
%$KS$ are not
%$N_G(S)$-cr.
\end{rems}

\section{An Important Example}
\label{sec:ex}

We consider a collection of important examples, which serve to illustrate
many of the points raised in the previous sections.
Throughout this section, we suppose that $p=2$ and
let $G$ be a simple group of type $G_2$.
We fix a maximal torus $T$ and a Borel subgroup $B$ of $G$ with $T\subseteq B$.
Let $\Psi$ be the set of roots of $G$ with respect to $T$.
We fix a base $\Sigma = \{\alpha,\beta\}$ for the set $\Psi^+$ of
positive roots with respect to $B$, where $\alpha$ is short and
$\beta$ is long.  The positive roots are $\alpha$, $\beta$,
$\alpha+\beta$, $2\alpha+\beta$, $3\alpha+\beta$ and $3\alpha+2\beta$.
For each root $\gamma$,
%we set $U_\gamma$ to be the
%corresponding root group and
%$\uu_\gamma$ to be $\Lie U_\gamma$;
we choose an isomorphism $\kappa_\gamma\colon k\ra U_\gamma$ and set
$s_\gamma=\kappa_\gamma(1)\kappa_{-\gamma}(-1)\kappa_\gamma(1)$ (cf.\
\cite[32.3]{Hum}).  Then $s_\gamma$ represents the reflection
corresponding to $\gamma$ in the Weyl group $N_G(T)/T$ of $G$.  Since
$p=2$, the order of $s_\gamma$ is 2 for every $\gamma\in \Psi$.

We use various equations from \cite[33.5]{Hum}.  Some are reproduced
below.  For brevity, we do not give the commutation relations between
the root subgroups: these are the equations of the form
$$ \kappa_\gamma(a)\kappa_{\gamma'}(b) =
\kappa_{\gamma'}(b)\kappa_\gamma(a)g,
$$ where $g$ is a product of elements of the form
$\kappa_{\gamma''}(p_{\gamma''}(a,b))$ over certain roots $\gamma''$,
each $p_{\gamma''}$ being a monomial in $a$ and $b$.  (Recall,
however, that $U_\gamma$ and $U_{\gamma'}$ commute if no positive
integral combination of $\gamma$ and $\gamma'$ is a root.)  We refer
to these equations collectively below as ``the CRs''.

%For $\gamma\in \Psi$, we set
%$$
%G_\gamma= \langle U_\gamma\cup U_{-\gamma}\rangle,\ \gg_\gamma
%= \Lie G_\gamma.
%$$
Since $G$ is simply connected, we have $G_\gamma\cong \SL_2(k)$ for
every $\gamma\in \Psi$, by Lemma~\ref{lem:sc}.

We have
\begin{equation}\label{eqn:fundalphaveepairings}
 \langle \alpha,\alpha^\vee\rangle=2,\ \langle
 \beta,\alpha^\vee\rangle=-3,\ \langle \alpha,\beta^\vee\rangle=-1,\
 \langle \beta,\beta^\vee\rangle=2
\end{equation}
(see \cite[32.3]{Hum}; note that $\inprod{\beta,\alpha}$ in
Humphreys's notation coincides with $\inprod{\beta,\alpha^\vee}$).
This yields
\begin{equation}
\label{eqn:alphaveepairings}
 \langle \alpha+\beta,\alpha^\vee\rangle=-1,\ \langle
 2\alpha+\beta,\alpha^\vee\rangle=1,\ \langle
 3\alpha+\beta,\alpha^\vee\rangle=3,\ \langle
 3\alpha+2\beta,\alpha^\vee\rangle=0.
\end{equation}
We have
$$ s_\alpha\cdot \alpha= -\alpha ,\ s_\alpha\cdot \beta=
3\alpha+\beta,
$$ and it follows that

\begin{equation}
\label{eqn:corootaction}
 s_\alpha\cdot \alpha^\vee= -\alpha^\vee,\ s_\alpha\cdot \beta^\vee =
\alpha^\vee+\beta^\vee.
\end{equation}

We need to know how $s_\alpha$ acts on the $U_\gamma$.  We can choose
the homomorphisms $\kappa_\gamma$ so that $s_\alpha$ maps each
$U_\gamma$ to $U_{s_\alpha\cdot \gamma}$ by conjugation in the
following way (see \cite[33.1 and 33.5]{Hum}):
\begin{equation}
\label{eqn:srootgps}
 s_\alpha\kappa_\beta(a)s_\alpha = \kappa_{3\alpha+\beta}(a),\
 s_\alpha\kappa_{\alpha+\beta}(a)s_\alpha =
 \kappa_{2\alpha+\beta}(a),\ s_\alpha\kappa_{2\alpha+\beta}(a)s_\alpha
 = \kappa_{\alpha+\beta}(a),
\end{equation}
$$ \ s_\alpha\kappa_{3\alpha+\beta}(a)s_\alpha=\kappa_\beta(a), \
s_\alpha\kappa_{3\alpha+2\beta}(a)s_\alpha=\kappa_{3\alpha+2\beta}(a),
$$ and we can choose $0\neq e_\gamma\in \uu_\gamma$ for each positive
root $\gamma$ such that the adjoint action of $s_\alpha$ on the
$e_\gamma$ is given by
\begin{equation}
\label{eqn:sonLie}
 \Ad {s_\alpha}(e_\beta) = e_{3\alpha+\beta},\ \Ad
 {s_\alpha}(e_{\alpha+\beta}) = e_{2\alpha+\beta},\ \Ad
 {s_\alpha}(e_{2\alpha+\beta}) = e_{\alpha+\beta},
\end{equation}
$$ \ \Ad {s_\alpha}(e_{3\alpha+\beta}) = e_\beta,\ \Ad
 {s_\alpha}(e_{3\alpha+2\beta}) = e_{3\alpha+2\beta}.
$$ Using Eqn.\ \eqref{eqn:srootgps} and the CRs, we get
\begin{eqnarray}
\label{eqn:scommrelns}
 & & s_\alpha \kappa_\beta(a) \kappa_{\alpha+\beta}(a')
              \kappa_{2\alpha+\beta}(b) \kappa_{3\alpha+\beta}(b')
              \kappa_{3\alpha+2\beta}(c) s_\alpha \\ & = &
              \kappa_\beta(b') \kappa_{\alpha+\beta}(b)
              \kappa_{2\alpha+\beta}(a') \kappa_{3\alpha+\beta}(a)
              \kappa_{3\alpha+2\beta}(ab'+a'b+c) \nonumber.
\end{eqnarray}

Let
\[
L:=\langle G_\alpha\cup T\rangle.
\]
Let $P$ be the parabolic subgroup of $G$ that contains $B$ and has $L$
as a Levi subgroup.  The roots of $R_u(P)$ with respect to $T$ are
$\beta$, $\alpha+\beta$, $2\alpha+\beta$, $3\alpha+\beta$,
$3\alpha+2\beta$.

Let $S$ be the torus $\alpha^\vee(k^*)$; we have $S=G_\alpha\cap T$.
Fix $t\in S$ such that $t$ has order 3.  By Eqns.\
\eqref{eqn:fundalphaveepairings} and \eqref{eqn:alphaveepairings},
\begin{equation}
\label{eqn:taction}
 \mbox{$t$ acts trivially on $U_\beta$, $U_{3\alpha+\beta}$,
                             $U_{3\alpha+2\beta}$, $\uu_\beta$,
                             $\uu_{3\alpha+\beta}$,
                             $\uu_{3\alpha+2\beta}$, and}
\end{equation}
$$ \mbox{$t$ acts non-trivially on $U_\alpha$, $U_{\alpha+\beta}$,
$U_{2\alpha+\beta}$, $\uu_\alpha$, $\uu_{\alpha+\beta}$,
$\uu_{2\alpha+\beta}$.}
$$

Set
\[
H := \langle \{s_\alpha,t\} \rangle\subseteq G_\alpha.
\]
Note that $H \cong S_3$.  Since $G_\alpha\cong \SL_2(k)$ (Lemma
\ref{lem:sc}), $\alpha^\vee$ is an isomorphism from $k^*$ onto
$G_\alpha\cap T$.  Set
$$ z := d\alpha^\vee(1)\in \gg_\alpha = \Lie G_\alpha
$$ (where we regard $1$ as an element of the tangent space
$T_1(k^*)\cong k$); then $z \neq 0$ and $k\cdot z$ is the centre of
$\gg_\alpha$.  In particular, $H$ centralizes $z$.  If $\gamma\in
\Psi$, then
\begin{equation}
\label{eqn:hcent}
 \mbox{$U_\gamma$ centralizes $z$}\iff \mbox{2 divides $\langle
\gamma,\alpha^\vee\rangle$}.
\end{equation}

By the CRs, $G_\alpha$ commutes with $G_{3\alpha+ 2\beta}$.  The
multiplication map $G_\alpha \times G_{3\alpha+2\beta}\ra G_\alpha
G_{3\alpha+2\beta}$ is bijective because $G_\alpha \cap
G_{3\alpha+2\beta}$, being a proper normal subgroup of $G_\alpha$,
must be trivial, but this map is not an isomorphism.  For
$\zz(\gg_{3\alpha+2\beta})\subseteq \Lie(G_{3\alpha+ 2\beta}\cap T)$,
so we have $\zz(\gg_{3\alpha+2\beta})\subseteq \cc_\frakt(G_\alpha)
\subseteq \cc_\frakt(H)$, which equals $k\cdot z$, so $\gg_\alpha\cap
\gg_{3\alpha+2\beta}$ is non-empty.  In fact $\gg_\alpha\cap
\gg_{3\alpha+2\beta}= k\cdot z$, because $k\cdot z$ is the unique
proper non-zero normal subalgebra of $\Lie G_\alpha= \Sl_2(k)$.

Now set
\[
M := G_\alpha G_{3\alpha+ 2\beta}.
\]
Note that $T\subseteq M$ and $M$ is a semisimple maximal rank subgroup of $G$
of type $\widetilde A_1 A_1$.
Consequently, $M$ is $G$-ir.
Moreover, we have
\[
C_G(M)= \ker(\alpha)\cap \ker(3\alpha+2\beta)
      = \ker(\alpha)\cap \ker(2\beta)
      = \ker(\alpha)\cap \ker(\beta)
      = Z(G)
      = \{1\}.
\]
%As $M$ is regular and reductive, $M$ is $G$-cr \cite[Prop.\ 3.20]{BMR}.
%Since $C_G(M) = \{1\}$, $M$ is therefore $G$-ir.
It is easy to see that $\zz(\mm)=k\cdot z$.
%It follows also from the preceding paragraph that
%$$   \zz(\mm)= k\cdot z. $$
%[Need to define $\mm$: $\mm= \Lie M$.]

\begin{lem}
\label{lem:M}
With the notation as above, we have
 \begin{itemize}
  \item[(a)] $(G,M)$ is a reductive pair;
  \item[(b)] $N_G(M)=M$;
  \item[(c)] $M= C_G(z)$.
 \end{itemize}
\begin{comment}
  (a) $(G,M)$ is a reductive pair.\smallskip\\ (b)
  $N_G(M)=M$.\smallskip\\ (c) $M= C_G(z)$.
\end{comment}
\end{lem}

\begin{proof}
\begin{comment}
 \begin{itemize}
  \item[(a)] Let $V$ be the subspace of $\gg$ spanned by $\uu_{\pm
  \beta}$, $\uu_{\pm (\alpha+\beta)}$, $\uu_{\pm (2\alpha+\beta)}$ and
  $\uu_{\pm (3\alpha+\beta)}$.  Suppose $\gamma\in \Psi(M)$,
  $\gamma'\in \Psi$ and $b$ is a positive integer such that
  $\gamma'+b\gamma\in \Psi(M)$.  It is easy to check that $\gamma'\in
  \Psi(M)$.  It follows from \cite[Prop.\ 27.2]{Hum} that $V$ is an
  $\Ad(M)$-stable complement to $\mm$ in $\gg$, so $(G,M)$ is a
  reductive pair.
  \item[(b)] Let $g\in N_G(M)$.  Without loss of generality, we can
assume that $g$ normalizes $T$, so $g$ permutes $\Psi(M)=\{\pm \alpha,
\pm (3\alpha+2\beta)\}$.  As $\alpha$ is short and $3\alpha+2\beta$ is
long, $g$ must normalize $G_\alpha$ and $G_{3\alpha+2\beta}$.  But
$G_\alpha$ and $G_{3\alpha+2\beta}$ are rank one groups, so they have
no outer automorphisms.  Since $C_G(M)=\{1\}$, we therefore have that
$g\in M$, as required.
  \item[(c)] Since $C_G(z)^0\supseteq M$ and $M$ is $G$-ir, $C_G(z)^0$
is $G$-ir, so $C_G(z)^0$ is reductive.  This implies that $C_G(z)^0$
is generated by the root subgroups that it contains together with $T$.
It follows from Eqns.\ \eqref{eqn:fundalphaveepairings},
\eqref{eqn:alphaveepairings} and \eqref{eqn:hcent} that $C_G(z)^0= M$.
Hence $C_G(z)\subseteq N_G(C_G(z)^0)= N_G(M)= M$.  Thus $C_G(z) = M$.
 \end{itemize}
\end{comment}
  (a). Since $\Psi(M)$ is a closed subsystem of $\Psi$, this follows immediately from Lemma~\ref{lem:redpair}.
  \begin{comment}
  Let $V$ be the subspace of $\gg$ spanned by $\uu_{\pm \beta}$,
  $\uu_{\pm (\alpha+\beta)}$, $\uu_{\pm (2\alpha+\beta)}$ and
  $\uu_{\pm (3\alpha+\beta)}$.  Suppose $\gamma\in \Psi(M)$,
  $\gamma'\in \Psi$ and $b$ is a positive integer such that
  $\gamma'+b\gamma\in \Psi(M)$.  It is easy to check that $\gamma'\in
  \Psi(M)$.
Since $M$ is generated by $U_{\pm \alpha}$ and $U_{\pm (3\alpha+ 2\beta)}$, it follows from \cite[Prop.\ 27.2]{Hum} that $V$ is an
  $\Ad(M)$-stable complement to $\mm$ in $\gg$, so $(G,M)$ is a
  reductive pair.
  \end{comment}

  (b). Let $g\in N_G(M)$.
Without loss of generality, we can assume that $g$ normalizes $T$,
so $g$ permutes $\Psi(M)=\{\pm \alpha, \pm (3\alpha+2\beta)\}$.
As $\alpha$ is short and $3\alpha+2\beta$ is long,
$g$ must normalize $G_\alpha$ and $G_{3\alpha+2\beta}$.
But $G_\alpha$ and $G_{3\alpha+2\beta}$ are rank one groups,
so they have no outer automorphisms. %\smallskip\\
Since $C_G(M)=\{1\}$, we therefore have that $g\in M$, as required.

  (c). Since $C_G(z)^0\supseteq M$ and $M$ is $G$-ir, $C_G(z)^0$ is $G$-ir,
so $C_G(z)^0$ is reductive.
This implies that $C_G(z)^0$ is generated by the root
subgroups that it contains together with $T$.
It follows from Eqns.\ \eqref{eqn:fundalphaveepairings},
\eqref{eqn:alphaveepairings} and
\eqref{eqn:hcent} that $C_G(z)^0= M$.
Hence $C_G(z)\subseteq N_G(C_G(z)^0)= N_G(M)= M$, by part (b).
Thus $C_G(z) = M$.
\end{proof}

\begin{lem}
\label{lem:HGcr}
With the notation as above, we have
\begin{itemize}
\item[(a)] $H$ is $G$-completely reducible and $M$-completely reducible;
\item[(b)] $C_G(H)= G_{3\alpha+ 2\beta}$ and $N_G(H)=HG_{3\alpha+ 2\beta}$.
\end{itemize}
\begin{comment}
 (a) $H$ is $G$-completely reducible and $M$-completely reducible.\smallskip\\
 (b) $C_G(H)= G_{3\alpha+ 2\beta}$ and $N_G(H)=HG_{3\alpha+ 2\beta}$.
\end{comment}
\end{lem}

\begin{proof}
\begin{comment}
 \begin{itemize}
  \item[(a)] It is easily checked that $H$ is not contained in any
  Borel subgroup of $L$, so $H$ is $L$-ir.  Now $L$ is a Levi subgroup
  both of $G$ and of $M$, so $H$ is both $G$-cr and $M$-cr by
  \cite[Cor.\ 3.22]{BMR}.
 \item[(b)] Since $H$ is $G$-cr, $C_G(H)$ is $G$-cr by \cite[Cor.\
3.17]{BMR}, so $C_G(H)^0$ is reductive.  Now $C_G(H)^0$ cannot have
rank $2$, because $H$ is not centralized by any maximal torus of $G$.
Hence $C_G(H)^0$ must be equal to its rank $1$ subgroup $G_{3\alpha+
2\beta}$.  By a similar argument, $C_G(G_{3\alpha+ 2\beta})^0=
G_\alpha$.

By \cite[Thm. 3.10]{BMR}, $G_{3\alpha+ 2\beta}=C_G(H)^0$ is $G$-cr.
Let $g\in N_G(H)$.  Then $g$ normalizes $C_G(H)^0= G_{3\alpha+
2\beta}$, so $g$ normalizes $C_G(G_{3\alpha+ 2\beta})^0= G_\alpha$, so
$g$ normalizes $M$.  Lemma \ref{lem:M}(a) now implies that $N_G(H)=
N_M(H)= HG_{3\alpha+ 2\beta}$ and $C_G(H)= C_M(H)= G_{3\alpha+
2\beta}$.
 \end{itemize}
\end{comment}
  (a). It is easily checked that $H$ is not contained in any Borel
 subgroup of $L$, so $H$ is $L$-ir.  Now $L$ is a Levi subgroup both
 of $G$ and of $M$, so $H$ is both $G$-cr and $M$-cr by \cite[Cor.\
 3.22]{BMR}.%\smallskip\\

  (b). Since $H$ is $G$-cr, $C_G(H)$ is $G$-cr by \cite[Cor.\ 3.17]{BMR},
  so $C_G(H)^0$ is reductive.  Now $C_G(H)^0$ cannot have
 rank $2$, because $H$ is not centralized by any maximal torus of $G$.
 Hence $C_G(H)^0$ must be equal to its rank $1$ subgroup $G_{3\alpha+
 2\beta}$.  Now $G_{3\alpha+ 2\beta}=C_G(H)^0$ is $G$-cr by \cite[Thm.\ 3.10]{BMR},
 so $C_G(G_{3\alpha+ 2\beta})^0= G_\alpha$ by a similar argument to that for $C_G(H)^0$.

%By \cite[Thm.\ 3.10]{BMR}, $G_{3\alpha+ 2\beta}=C_G(H)^0$ is $G$-cr.
Let $g\in N_G(H)$.  Then $g$ normalizes $C_G(H)^0= G_{3\alpha+
2\beta}$, so $g$ normalizes $C_G(G_{3\alpha+ 2\beta})^0= G_\alpha$, so
$g$ normalizes $M$.  Lemma \ref{lem:M}(b) now implies that $N_G(H)=
N_M(H)= HG_{3\alpha+ 2\beta}$ and $C_G(H)= C_M(H)= G_{3\alpha+
2\beta}$.
\end{proof}

The following example shows that the converse of Corollary
\ref{cor:LGsep} can fail, even when $H$ is $L$-ir.  It also gives a
counterexample to the converse of the first assertion of Theorem
\ref{thm:redpairs} (note that $(G,L)$ is a reductive pair by Lemma \ref{lem:redpair}).

\begin{prop}
\label{prop:LsepGnonsep}
 The subgroup $H$ is separable in $L$, but not in $G$.
\end{prop}

\begin{proof}
Since $DL = G_\alpha\cong \SL_2(k)$,
$L$ is isomorphic either to $\GL_2(k)$ or to $\SL_2(k)\times k^*$.
To rule out the latter case, it's enough to show
that $s_\alpha$ acts non-trivially on $\frakt$;
this follows from Eqn.\ (\ref{eqn:corootaction}).
It now follows easily that $\cc_\frakl(H)= k\cdot z$.
Now $z$ is tangent to $\Lie Z(L)$ as $L\cong \GL_2(k)$.
We deduce that $H$ is separable in $L$.

It follows from Eqns.\ \eqref{eqn:sonLie} and \eqref{eqn:taction}
that $H$ centralizes
$e_\beta+ e_{3\alpha+ \beta}$.
But $C_G(H)= G_{3\alpha+ 2\beta}$ by Lemma~\ref{lem:HGcr}(b),
so $e_\beta+ e_{3\alpha+ \beta}$
is not tangent to $\Lie C_G(H)= \gg_{3\alpha+ 2\beta}$.
Hence $H$ is not separable in $G$.
\end{proof}

\begin{rem}
 It is easily checked that for every semisimple $x\in \cc_\gg(H)$, we have
 $x\in \Lie C_G(H)$.  The same result cannot hold if we replace $H$ with a $G$-irreducible and
 non-separable subgroup of $G$; cf.\ the proof of
 \cite[Thm.~3.39]{BMR}.
\end{rem}

%\begin{rem}
%Note that although $H$ is separable in $L$,
%$\pi_L(H)$ is not separable in $L_{\rm ad}\cong \PGL_2(k)$.
%For $C_{L_{\rm ad}}(\pi_L(H)) = \{1\}$,
%but $\pi_L(H)$ centralizes a non-zero
%semisimple element of $\Lie L_{\rm ad}$.
%\end{rem}

For $a\in k$, set
$$ u(a)= \kappa_\beta(a)\kappa_{3\alpha+\beta}(a).
$$ The CRs yield
\begin{equation}
\label{eqn:umult}
 u(a)u(b)= u(a+b)\kappa_{3\alpha+2\beta}(ab)
\end{equation}
and
\begin{equation}
\label{eqn:uinv}
 u(a)^{-1}= u(a)\kappa_{3\alpha+2\beta}(a^2).
\end{equation}
Define
\[
H_a= u(a)Hu(a)^{-1}.
\]
By Eqn.\ \eqref{eqn:taction}, $u(a)$ centralizes $t$; by Eqn.\
\eqref{eqn:scommrelns}, we have $u(a)s_\alpha u(a)^{-1}= s_\alpha
\kappa_{3\alpha+2\beta}(a^2)$.  Hence
$H_a = \langle \{ t, s_\alpha \kappa_{3\alpha+ 2\beta}(a^2) \}\rangle
\subseteq M$ for all $a\in k$.

The following example shows that Theorem \ref{thm:finorbit}(b) can
fail if we do not require $K$ to be separable in $G$ (recall that
$(G,M)$ is a reductive pair by Lemma \ref{lem:M}(a)).
This example also proves Theorem \ref{thm:countereg}.

\begin{exmp}
\label{exmp:orbitcountereg}
Let $(m_1,m_2)= (s_\alpha, t \kappa_{3\alpha+2\beta}(1))$.  For
$a\in k$, set
\[
(m_1(a),m_2(a))= u(a)\cdot (m_1,m_2).
\]
By construction,
$(m_1(a),m_2(a))$ is $G$-conjugate to $(m_1(b),m_2(b))$
for all $a,b\in k$.  We now show that if $a\neq b$, then
$(m_1(a),m_2(a))$ is not $M$-conjugate to
$(m_1(b),m_2(b))$.  This shows that $G\cdot (m_1,m_2)$ is
an infinite union of $M$-conjugacy classes.

For $a \in k$ let
\[
\widehat{H}_a := \langle \{m_1(a),m_2(a)\} \rangle.
\]
Note that $t$ commutes with $\kappa_{3\alpha+2\beta}(1)$ by Eqn.\ (\ref{eqn:taction}), so
\[
\widehat{H}_0 =  \langle \{ m_1, m_2 \} \rangle= \langle \{ s_\alpha, t, \kappa_{3\alpha+2\beta}(1) \} \rangle = \langle H\cup \{\kappa_{3\alpha+2\beta}(1)\}\rangle
\]
and $\widehat{H}_a = \langle H_a\cup \{\kappa_{3\alpha+2\beta}(1)\}\rangle$ similarly.  Hence
\[
C_G(\widehat{H}_0) = C_G(H)\cap C_G(\kappa_{3\alpha+2\beta}(1)) =
           G_{3\alpha+2\beta}\cap C_G(\kappa_{3\alpha+2\beta}(1)) =
           U_{3\alpha+2\beta},
\]
by Lemma \ref{lem:HGcr}(b), so
\[
C_G(\widehat{H}_a)= C_G(u(a)\widehat{H}_0u(a)^{-1})= u(a)C_G(\widehat{H}_0)u(a)^{-1}= u(a)U_{3\alpha+2\beta}u(a)^{-1}=
U_{3\alpha+2\beta}.
\]
Now let $a, b \in k$, and suppose that $(m_1(a),m_2(a))$ and
$(m_1(b),m_2(b))$ are $M$-conjugate.
Then there exists $m\in M$
such that $(mu(a))\cdot (m_1,m_2)=u(b)\cdot (m_1,m_2)$.  We
have $mu(a)u(b)^{-1}\in C_G(\langle \{m_1(b),m_2(b)\}\rangle)=
C_G(\widehat{H}_b)= U_{3\alpha+2\beta}\subseteq M$, so
$u(a)u(b)^{-1}\in M$.  But $u(a)u(b)^{-1}= u(a+b)
\kappa_{3\alpha+2\beta}(ab+b^2)$ by Eqns.\ (\ref{eqn:umult}) and
(\ref{eqn:uinv}), so $u(a+b)\in M$.  So we must have $u(a+b)=1$,
whence $a=b$.
Thus if $a \neq b$, then $(m_1(a),m_2(a))$ and
$(m_1(b),m_2(b))$ are not $M$-conjugate.

These calculations show that even though $(G,M)$ is a reductive pair,
$G\cdot(m_1,m_2) \cap M^2$ consists of an infinite union of $M$-conjugacy classes.
Observe that this is consistent with Theorem \ref{thm:finorbit}(b); a
similar calculation to the one in the proof of Proposition
\ref{prop:LsepGnonsep} shows that $\widehat{H}_0$
is not separable in $G$.
%This example proves Theorem \ref{thm:countereg}.
\end{exmp}

\begin{rem}
 The $n$-tuple $(m_1,m_2,\ldots, m_n)$ yields a similar example for any $n\geq 2$, where $m_1$ and $m_2$ are as above and $m_3= \cdots = m_n=1$.
\end{rem}

\begin{comment}
K\"ulshammer [ref] asked the following question.

\begin{question}
 Let $G$ be an algebraic group and let $F$ be a finite group.  Fix a
Sylow $p$-subgroup $F_p$ of $F$.  Given a homomorphism $\rho_0\colon
F\ra G$, is it true that there are only finitely many $G$-conjugacy
classes of homomorphisms $\rho\colon F\ra G$ such that the
restrictions $\rho|_{F_p}$ and $\rho_0|_{F_p}$ are $G$-conjugate?
\end{question}

The answer is yes in many cases: if the characteristic is zero or if
$G=\GL_n(k)$, for example [refs].  The answer is no for general
non-reductive groups [ref].
\end{comment}

%We now give a counterexample of the failure of the forward implication
%of Proposition \ref{prop:regGcr} in bad characteristic (see Subsection
%\ref{subsec:reg}).
Our next example shows that the second assertion
of Theorem~\ref{thm:redpairs} is false with the separability
assumption on $K$ removed, even though $(G,M)$ is a reductive pair by
Lemma \ref{lem:M}(a).
Since $M$ is a regular subgroup of $G$, this is also a new example of the failure of
\cite[Thm.\ 3.26]{BMR} in bad characteristic (compare
\cite[Ex.\ 3.45]{BMR}, which gives subgroups $H' \subseteq M' \subseteq G'$,
with $G'$ connected reductive, $M'$ regular in $G'$ such that
$H'$ is $M'$-cr, but not $G'$-cr).

\begin{prop}
\label{prop:GcrnotMcr}
Let $a\in k^*$.  Then $H_a$ is $G$-completely reducible but not
$M$-completely reducible.
\end{prop}

\begin{proof}
 Since $H_a$ is $G$-conjugate to $H$ and $H$ is $G$-cr (Lemma
\ref{lem:HGcr}(a)), $H_a$ is $G$-cr.
%To show that
%$H_a$ is not $M$-cr, we use the results and notation of
%\cite[Sec.~2.3]{BMR}.
Let $\lambda= \alpha^\vee+ 2\beta^\vee\in Y(T)$: then $\langle
\alpha,\lambda\rangle= 0$ and $\langle \beta,\lambda\rangle= 1$.  It
is clear that $P=P_\lambda$ and $L = L_\lambda$.  We have a
homomorphism $c_\lambda\colon P_\lambda\ra L_\lambda$ as defined in
Subsection \ref{subsec:noncon}.  If $h\in H$ then $u(a)hu(a)^{-1}= hu$ for
some $u\in U_{3\alpha+ 2\beta}\subseteq \ker(c_\lambda)$, so
$c_\lambda(H_a)= H$.  To prove that $H_a$ is not $M$-cr, it suffices
by \cite[Lem.\ 2.17 and Thm.\ 3.1]{BMR} to show that $H_a$ and $H$ are
not $M$-conjugate.  But this is the case since $u(a)\not\in M$ and
$N_G(H)\subseteq M$ (Lemma \ref{lem:HGcr}(b)), so we are done.
\end{proof}

\begin{rem}
\label{rem:thm1.3(c)fails}
Proposition \ref{prop:GcrnotMcr} shows that
part (c) of Theorem \ref{thm:finorbit} fails
without the separability hypothesis on $K$:
for, by Theorem \ref{thm:gcrcrit},
the $G$-conjugacy class $G\cdot (m_1(a),m_2(a))$ is closed
in $G \times G$ but $M\cdot (m_1(a),m_2(a))$ is not.
\end{rem}

%\begin{rem}
%The subgroup $H_0 = H$ is $L$-ir, so it is $M$-cr by
%\cite[Cor.\ 3.22]{BMR}.
%$L$ is also a Levi of $M$.
%\end{rem}

\begin{comment}
% Following result is probably true, but isn't useful.
\begin{lem}
\label{lem:Hmaximality}
Let $K$ be a proper connected regular reductive subgroup of $G$ such
that $H\subseteq G$.  Then some conjugate of $K$ is contained in
$G_\alpha G_{3\alpha+2\beta}$.  In particular, every connected
unipotent subgroup of $K$ is abelian.
\end{lem}

\begin{proof}
 The semisimple element $t$ lies in a maximal torus of $K$.  Since $K$
is regular, this torus is also a maximal torus of $G$, so we can
assume without loss of generality that this torus is $T$.

 Any root is $W$-conjugate either to $\alpha$ or $3\alpha+2\beta$, so
if $K$ has semisimple rank at most $1$, then $K$ is conjugate to one
of $G_\alpha$, $G_{3\alpha+2\beta}$, $\langle G_\alpha$, $T\rangle$,
$\langle G_{3\alpha+2\beta},T\rangle$ or to a subtorus of $T$.  So
assume that $K$ has semisimple rank 2 and that $T\subseteq K$.  ...
\end{proof}
\end{comment}

The next example shows that Proposition \ref{prop:sepovergp} can fail
if we allow $H$ to be non-$G$-cr.  First we need a refinement of Lemma \ref{lem:sep}.
If $G_1$ is a reductive group and $H_1$ is a subgroup of $G_1$, then we say
that $x\in \gg_1$ is a {\em witness to the non-separability of $H_1$} if
$x\in \cc_{\gg_1}(H_1)$ but $x\not\in \Lie C_{G_1}(H_1)$.

\begin{lem}
\label{lem:nilptwitness}
 Let $f\colon G_1\ra G_2$ be an epimorphism of reductive groups such
that $\ker df$ consists of semisimple elements and let $H_1$ be a
subgroup of $G_1$.  Let $x\in \cc_{\gg_1}(H_1)$ be nilpotent.
Then $x$ is a witness to the non-separability of $H_1$ if and only
if $df(x)$ is a witness to the non-separability of $f(H_1)$.
\end{lem}

\begin{proof}
 It is clear that $df(x)\in \cc_{\gg_2}(f(H_1))$ and that
$df(x)$ is tangent to $C_{G_2}(f(H_1))$ if $x$ is tangent
to $C_{G_1}(H_1)$.  Conversely, suppose that $df(x)$ is
tangent to $C_{G_2}(f(H_1))$.  By \cite[Prop.~IV.14.26]{borel},
there exists a connected unipotent subgroup $U_2$ of $C_{G_2}(f(H_1))$
such that $df(x)\in \Lie U_2$.  By \cite[V.22.1]{borel} the restriction
of $f$ is an isogeny from $D(G_1^0)$ to $D(G_2^0)$.  Let $U_1= (f^{-1}(U_2)\cap D(G_1^0))^0$.
Any semisimple element of $U_1$ must belong to the finite group $\ker f\cap D(G_1^0)$, so $U_1$
has only finitely many semisimple elements.  Hence $U_1$ is a unipotent group.  As $\ker df$
consists of semisimple elements and $U_2\subseteq D(G_2^0)$, it follows that the restriction
of $f$ is an isomorphism from $U_1$ onto $U_2$.  Hence we can choose $x'\in \Lie U_1$ such that
$df(x')=df(x)$.  Now $\ker df\subseteq \zz(\gg_1)$, by \cite[V.22.2]{borel}, and $x,x'$ are both
nilpotent, so we must have $x=x'$, whence $x\in \Lie U_1$.

To complete the proof, it suffices to show that $U_1\subseteq C_{G_1}(H_1)$.
Fix $h\in H_1$.  If $u\in U_1$, then $f(h)f(u)f(h)^{-1}= f(u)$,
so we have $huh^{-1}=cu$ for some $c\in \ker f$.
The map $u\mapsto huh^{-1}u^{-1}$
is therefore a morphism from the connected
set $U$ to the finite set $\ker f\cap D(G_1^0)$,
so $huh^{-1}u^{-1}=1$ for all $h\in H_1$ and all $u\in U_1$.
Hence $U_1\subseteq C_{G_1}(H_1)$, as required.
\end{proof}

\begin{exmp}
\label{exmp:nonGcr}
 Let $C= \{u(a) \mid a\in k\}$.  Let $H'= \langle H\cup C\rangle$ and
suppose that $K$ is any reductive subgroup of $G$ containing $H'$.
Set $y:=e_\beta+ e_{3\alpha+ \beta}$.  Then $y$ belongs to the tangent
space $T_1(C)$, which is contained in $\Lie K$.  Now $\uu_{3\alpha+
2\beta}$ is centralized by $U_\beta$ and $U_{3\alpha+ 2\beta}$ and we
have $\Ad \kappa_{3\alpha+\beta}(a)(e_\beta) = e_\beta+ ae_{3\alpha+
2\beta}$, $\Ad \kappa_\beta(a)(e_{3\alpha+ \beta}) = e_{3\alpha+
\beta}+ ae_{3\alpha+2\beta}$ by the CRs, so $\Ad u(a)(y)= y$.  Hence $y\in \cc_\kk(H')$, but $y$ is not
tangent to $C_K(H')$ since it is not tangent to $C_G(H)$ (cf.\ the
proof of Proposition \ref{prop:LsepGnonsep}).  By Lemma \ref{lem:nilptwitness}, $y$ is a witness to the non-separability of $\pi_K(H')$, so $\pi_K(H')$ is not separable in $K_{\rm ad}$.
\end{exmp}

Here is a further example arising from this construction which relates
to the discussions in Sections \ref{sec:redpairs} and \ref{sec:cent+norm}
(see in particular
Remarks \ref{rem:C_G(S)redpair}(ii) and \ref{rem:KvsKN}(ii)).

\begin{exmp}
\label{exmp:KnotKS}
Recall that $S$ is the torus $\alpha^\vee(k^*)$
and therefore, $(G,N_G(S))$ is a reductive pair by
Proposition~\ref{prop:C_G(S)redpair}(b).
We have that $H_a \subseteq N_G(S)$, since $U_{3\alpha+2\beta}$
centralizes $S$.  Let $a\in k^*$. We will show that $H_aS$ is not
$G$-cr, although $H_a$ is (Proposition \ref{prop:GcrnotMcr}).  Let
$\lambda\in Y(T)$ and $c_\lambda$ be as in the proof of Proposition
\ref{prop:GcrnotMcr}.  Then
$$ c_\lambda(H_aS)= \langle S\cup\{s_\alpha\} \rangle.
$$
If $H_aS$ lies in a Levi subgroup of $P_\lambda$ then $uH_aSu^{-1}\subseteq L_\lambda$ for some $u\in R_u(P_\lambda)$, so $c_\lambda(H_aS)= uH_aSu^{-1}$.  Thus to show that $H_aS$ is not $G$-cr, it is enough to show that
$\langle S\cup\{s_\alpha\} \rangle$ is not $R_u(P_\lambda)$-conjugate
to $H_aS$.  To see this, note that if $u\in R_u(P_\lambda)$ with
$u\langle S\cup\{s_\alpha\} \rangle u^{-1}= H_aS$, then $uSu^{-1}= S$,
whence $u$ centralizes $S$ (since $S$ normalizes $R_u(P_\lambda)$).
%(see \cite[Thm.\ 10.6(5)(ii)]{borel}). %apply it to $B$
Since the centralizer of $S$ in $R_u(P_\lambda)$ is $U_{3\alpha+
2\beta}$ by Eqns.\ (\ref{eqn:fundalphaveepairings}) and
(\ref{eqn:alphaveepairings}), we have $u\in U_{3\alpha+ 2\beta}$.  But
$U_{3\alpha+ 2\beta}$ centralizes $H_aS$, so $\langle
S\cup\{s_\alpha\} \rangle= H_aS$, a contradiction.

Since $H_aS$ is not $G$-cr, $H_aS$ is not $N_G(S)$-cr (Corollary \ref{cor:KS<=>K}).  Since the canonical projection $f\colon N_G(S)\ra N_G(S)/S$ is non-degenerate and $f(H_a)= f(H_aS)$, it follows from \cite[Lem.\ 2.12(ii)]{BMR} that $H_a$ is not $N_G(S)$-cr.  Thus we have an example of a subgroup
$H_a \subseteq N_G(S)$, with $S$ linearly reductive --- in fact, a
torus --- such that $H_a$ is $G$-cr, but not $N_G(S)$-cr.

% MB:- I'm not sure of the point of the following???
% This can also be seen slightly more directly
% set $D= \langle S,s_\alpha,G_{3\alpha+ 2\beta} \rangle$.
% Then, by equations \eqref{eqn:fundalphaveepairings} and \eqref{eqn:alphaveepairings},
% $D^0= \langle S,G_{3\alpha+ 2\beta} \rangle= C_G(S)^0$, and $D$ normalizes $S$.
% Thus $D \subseteq N_G(S)$ and $D$ has finite index in $N_G(S)$, so by
% \cite[Prop. 2.15]{BMR2}, a subgroup of $D$ is $D$-cr if and only if it is $N_G(S)$-cr.
% Now the image of $H_a$ under the canonical projection from
% $D$ to $D/\langle S,s_\alpha \rangle = G_{3\alpha+ 2\beta}$
% has order $2$ and so is not $G_{3\alpha+ 2\beta}$-cr.
\end{exmp}

%\section{Rationality}\label{sec:fields}

Finally, we consider a rationality question.  Let $k_0$ be a subfield
of an algebraically closed field $k_1$.  Suppose that $G_1$ is a
reductive algebraic group defined over $k_0$.  If $K_1$ is a subgroup
of $G_1$ defined over $k_0$, then we say that $K_1$ is
\emph{$G_1$-completely reducible over $k_0$} if whenever $P_1$ is an
R-parabolic subgroup of $G_1$ such that $K_1\subseteq P_1$ and $P_1$
is defined over $k_0$, there exists an R-Levi subgroup $L_1$ of $P_1$
such that $K_1\subseteq L_1$ and $L_1$ is defined over $k_0$ (see
\cite[Sec.~5]{BMR} for further details).  In particular, $K_1$ is
$G_1$-cr if and only if $K_1$ is $G_1$-cr over $k_1$.  An example of McNinch
\cite[Ex.\ 5.11]{BMR} shows that if $K_1$ is $G_1$-cr over $k_0$, then $K_1$
need not be $G_1$-cr over $k_1$.  The next example shows that the converse
can also happen.

\begin{comment}
\begin{lem}
\label{lem:levirat}
 Let $P$ be a parabolic subgroup of $G$, let $L$ be a Levi subgroup of
$P$ and suppose that $P,L$ are defined over $k_0$.  Given $p\in P$,
write $p=lu$, where $l\in L$ and $u\in R_u(P)$.  Then $p\in P(k_0)$ if
and only if $l\in L(k_0)$ and $u\in R_u(P)(k_0)$.
\end{lem}

\begin{proof}
 It is immediate that $p\in P(k_0)$ if $l\in L(k_0)$ and $u\in
R_u(P)(k_0)$.  Conversely, suppose that $p\in P(k_0)$.  By
\cite[Prop.\ 20.5, Thm.\ 6.8]{borel}, $R_u(P)$ is defined over $k_0$
and the canonical projection $\pi_P\colon P\ra P/R_u(P)$ is defined
over $k_0$.  Hence the map $\psi= \pi_P\circ i$ is defined over $k_0$,
where $i$ is the inclusion of $L$ in $P$.  So $l=\psi(p)\in L(k_0)$
and $u= l^{-1}p\in R_u(P)(k_0)$, as required.
\end{proof}
\end{comment}

\begin{exmp}
\label{exmp:mcninchcvse}
% Let $p=2$, $G$, $H$, etc., be as above. % in Section \ref
%{sec:ex}.
Suppose that $k_0$ is a subfield of $k$ such that $G$ is defined over
$k_0$ and $k_0$-split.  We can assume that $T$ is chosen so that $T$
is defined over $k_0$ and $k_0$-split and so that for every $\gamma\in
\Psi$, the homomorphisms $\gamma\colon T\ra k^*$, $\gamma^\vee\colon k^*\ra T$ and
$\kappa_\gamma\colon k\ra U_\gamma$ are defined over $k_0$.  Now
suppose that $k/k_0$ is not separable.  Then $k_0$ is not perfect, so
there exists $a\in k\setminus k_0$ such that $a^2\in k_0$.  Consider
$H_a= \langle \{t,s_\alpha \kappa_{3\alpha+ 2\beta}(a^2) \}\rangle$.  Then
$H_a$ is defined over $k_0$, since $H_a$ is a finite subgroup of $G(k_0)$.  As $H_a$ is
$G$-cr (Proposition \ref{prop:GcrnotMcr}), $H_a$ is $G$-cr over $k$.  We
show that $H_a$ is not $G$-cr over $k_0$.

\begin{comment}
Let $P$ be a parabolic subgroup of $G$ such that $P$ is defined over
$k_0'$ and $H\subseteq P$.  We may assume that $P\neq G$.  Since $H$
is $G$-cr, there is a Levi subgroup $L$ of $P$ such that $H\subseteq
P$.  Choose indivisible $\mu\in Y(G)$ such that $P= P_\mu$ and
$L=L_\mu$.  We can regard $\mu$ as an element of $Y(G_{3\alpha+
2\beta})$, and $\mu$ is non-trivial as $P\neq G$.  Since $G_{3\alpha+
2\beta}\cong \PGL_2(k)$, there exists $z\in G_{3\alpha+ 2\beta}$ such
that $\mu= z\cdot \lambda$.  Hence $zP_\lambda z^{-1}= P_\mu= P$.  By
\cite[Thm.\ 20.9(iii)]{borel}, there exists $g\in G(k_0')$ such that
$gP_\lambda g^{-1}= P$.  We can write $g= zp$ for some $p\in
N_G(P_\lambda)= P_\lambda$.

Let $P$ be a parabolic subgroup of $G$ such that $P$ is defined over
$k_0'$ and $H\subseteq P$.  We may assume that $P\neq G$.  Since any
Levi subgroup of a Borel subgroup of $G$ is abelian and $H$ is
non-abelian, $P$ is not a Borel subgroup of $G$.  By \cite[?]{borel},
$P\cap P_\lambda$ contains a maximal torus $S$ of $G$ such that $S$ is
defined over $k_0'$.  There exists $u\in R_u(P_\lambda)$ such that
$uSu^{-1}= T$.  Hence $uPu^{-1}$ is a parabolic subgroup of $G$
containing $T$ and $uHu^{-1}$.
\end{comment}

Recall that $P = P_\lambda$ and $L = L_\lambda$, where $\lambda$ is as
in the proof of Proposition \ref{prop:GcrnotMcr}.  Since $G$ and $T$ are split, $P$
and $L$ are defined over $k_0$ (cf.\ \cite[Props.~V.20.4 and V.20.5]{borel}).  Suppose there exists a Levi subgroup
$L'$ of $P$ such that $L'$ is defined over $k_0$ and $H_a$ is
contained in $L'$.  By \cite[Prop.~V.20.5]{borel}, there exists $u\in
R_u(P)(k_0)$ such that $L'= uLu^{-1}$.  Then $u^{-1}H_au$,
$u(a)^{-1}H_au(a)\subseteq L$, so $c_\lambda(u^{-1}gu)= u^{-1}gu$ and
$c_\lambda(u(a)^{-1}gu(a))= u(a)^{-1}gu(a)$ for all $g\in H_a$.  Since
$c_\lambda(u)= c_\lambda(u(a))= 1$, we have
\[
u^{-1}gu = c_\lambda(u^{-1}gu) = c_\lambda(g)
%= c_\lambda(u(a)(u(a)^{-1}gu(a))u(a)^{-1})
= c_\lambda(u(a)^{-1}gu(a)) =
u(a)^{-1}gu(a)
\]
for all $g\in H_a$.  Thus $u= u(a)c$ for some $c \in C_G(H)\cap
R_u(P)= G_{3\alpha+ 2\beta}\cap R_u(P) = U_{3\alpha+ 2\beta}$ (using Lemma \ref{lem:HGcr}(b)): say $u=
\kappa_\beta(a)\kappa_{3\alpha+ \beta}(a) \kappa_{3\alpha+ 2\beta}(y)$
for some $y\in k$.  But then $u\not\in P(k_0)$ since $a\not\in k_0$, a
contradiction.  Thus no such $L'$ can exist, and $H_a$ is not
$G$-cr over $k_0$.
\begin{comment}
We may regard $u$ as an element of $P_{\alpha^\vee}$, with $z\in
L_{\alpha^\vee}$ and $u(a)\in R_u(P_{\alpha^\vee})$.  By \cite[20.5
Prop.,\ 6.8 Thm.]{borel}, $R_u(P_{\alpha^\vee})$ is defined over $k_0$
and the canonical projection $\pi\colon P_{\alpha^\vee}\ra
P_{\alpha^\vee}/R_u(P_{\alpha^\vee})$ is defined over $k_0$.  Hence
the isomorphism $\psi= \pi\circ i$ from $L_{\alpha^\vee}$ to
$P_{\alpha^\vee}/R_u(P_{\alpha^\vee})$ is defined over $k_0$, where
$i$ is the inclusion of $L_{\alpha^\vee}$ in $P_{\alpha^\vee}$.  We
have $z= \psi^{-1}(\pi(u))\in L_{\alpha^\vee}(k_0)$, so $u\in
P_{\alpha^\vee}(k_0)$.  But $u\not\in P_{\alpha^\vee}(k_0)$, because
$a\not \in k_0$, a contradiction.  Thus no such $L$ can exist, and
$H_a$ is not $G(k_0)$-cr.
\end{comment}
\end{exmp}

\bigskip {\bf Acknowledgements}: The authors acknowledge the financial
support of EPSRC Grant EP/C542150/1 and Marsden Grant UOC0501.  Part
of the research for this paper was carried out while the authors were
staying at the Mathematical Research Institute Oberwolfach supported
by the ``Research in Pairs'' programme.  We are grateful to Martin
Liebeck and Gary Seitz for helpful discussions.
%We are grateful to G.\ Seitz for pointing out the case used in
%Example \ref{ex:e6}.

\bigskip

%%%%%%%%%%%%%%%%%%%%%%%%%%%%%%%%%%%%%%%%%%%%%%%%%%%%%%%%%%%%%%%%%%%%%%
%%%%%%%%%%%%% bibliography
%%%%%%%%%%%%%%%%%%%%%%%%%%%%%%%%%%%%%%%%%%%%%%%%%%%%%%%%%%%%%%%%%%%%%%

\end{document}